\documentclass[12pt]{article}
\usepackage{hyperref}
\usepackage{amsmath}
\usepackage{graphicx}
\usepackage{amssymb} 
 \usepackage{xcolor}
\definecolor{myblue}{rgb}{0.14,0.11,0.49}
\definecolor{myred}{rgb}{0.74,0.22,0.15}
\definecolor{mygreen}{rgb}{0.05,0.52,0.42}
\definecolor{myyellow}{rgb}{0.96,0.92,0.13}
\definecolor{myorange}{rgb}{1,0.61,0.36}
\definecolor{mypurple}{rgb}{0.71,0.02,1}
\definecolor{noir}{gray}{0.} 

\newcommand{\Mat}[1]{{{\boldsymbol{#1}}}}
\newcommand{\abs}[1]{\left\vert#1\right\vert}
\def\be{\begin{equation}}
\def\ee{\end{equation}}
\def\bea{\begin{eqnarray}}
\def\eea{\end{eqnarray}}
\def\bi{\begin{itemize}}
\def\ei{\end{itemize}}
\def\noi{\noindent}
\def\dd{\mathrm{d}}

\date{}

\title{Defining the space in a general spacetime}

\author{
Mayeul Arminjon$^{1,2}$\\
$^1$ \small\it  \ Univ. Grenoble Alpes, Lab. 3SR, F-38000 Grenoble, France. \\
$^2$ \small\it CNRS, Lab. 3SR, F-38000 Grenoble, France.\\
\small E-mail: Mayeul.Arminjon@3sr-grenoble.fr}

\begin{document}
\maketitle

\begin{abstract} 
\noi A global vector field $v$ on a ``spacetime" differentiable manifold $\mathrm{V}$, of dimension $N+1$, defines a congruence of world lines: the maximal integral curves of $v$, or orbits. The associated global space $\mathrm{N}_v$ is the set of these orbits. A ``$v$-adapted" chart on $\mathrm{V}$ is one for which the $\mathbb{R}^N$ vector ${\bf x}\equiv (x^j)\ (j=1,...,N)$ of the ``spatial" coordinates remains constant on any orbit $l$. We consider non-vanishing vector fields $v$ that have non-periodic orbits, each of which is a closed set. We prove transversality theorems relevant to such vector fields. Due to these results, it can be considered plausible that, for such a vector field, there exists in the neighborhood of any point $X\in \mathrm{V}$ a chart $\chi $ that is $v$-adapted and ``nice", i.e., such that the mapping $\bar{\chi }: l\mapsto {\bf x}$ is injective --- unless $v$ has some ``pathological" character. This leads us to define a notion of ``normal" vector field. For any such vector field, the mappings $\bar{\chi }$ build an atlas of charts, thus providing $\mathrm{N}_v$ with a canonical structure of differentiable manifold (when the topology defined on $\mathrm{N}_v$ is Hausdorff, for which we give a sufficient condition met in important physical situations). Previously, a local space manifold $\mathrm{M}_\mathrm{F}$ had been associated with any ``reference frame" $\mathrm{F}$, defined as an equivalence class of charts. We show that, if $\mathrm{F}$ is made of nice $v$-adapted charts, $\mathrm{M}_\mathrm{F}$ is naturally identified with an open subset of the global space manifold $\mathrm{N}_v$. \\

\noi {\it Keywords:} Physical space; global vector field; reference fluid; orbit space; adapted chart; differentiable manifold; Kruskal-Szekeres coordinates.

\end{abstract} 

\section{Introduction}\label{Intro}

\subsection{Physical motivation}

The theory of relativity says that space and time merge into ``a kind of union of the two" (in Minkowski's words): the spacetime. However, the notion of a {\it physical space} should be useful also in relativistic physics. In our opinion it is even needed, for the following two reasons. (i) In experimental/observational work, one of course needs to define the spatial position of the experimental apparatus and/or of the observed system, and this is true also if the relativistic effects have to be considered. (ii) In quantum mechanics, the ``state" of a quantum-mechanical particle is a function $\psi $ of the position $x$ belonging to some 3-D ``physical space" $\mathrm{M}$, and taking values in $\mathbb{C}$ (or in a complex vector bundle). Note that defining such a space as a spacelike 3-D {\it submanifold of the spacetime manifold} (e.g. \cite{Lachièze-Rey2003}) can work to define an initial condition for a field in space-time, but does not allow one to define a spatial position in the way that is needed in the two foregoing examples: in those, one needs to identify spatial points that exist at least for some open interval of time --- e.g. to state that some objects maintain a fixed spatial position in some reference frame, or to define the stationary states of a quantum particle. In practice, the spatial position is taken to be the triplet of the spatial coordinates, ${\bf x}\equiv (x^j) \,(j=1,2,3)$. However, {\it a priori,} ${\bf x}$ {\it does not have a precise geometric meaning in a theory starting from a spacetime structure.} Only a notion of ``spatial tensors" has been defined for a general spacetime of relativistic gravity. This definition was based on the concept of ``reference fluid" \cite{Cattaneo1958, Massa1974a,Massa1974b}, also named ``reference body" \cite{Mitskievich1996} --- i.e., a three-dimensional congruence of world lines, whose the tangent vector is assumed to be a time-like vector field. The latter can be normed to become a four-velocity field $v$, i.e. $\Mat{g}(v,v)=1$ where $\Mat{g}$ is the spacetime metric. The data of the four-velocity field $v$ 
allows one to define the spatial projection operator $\Pi _X$ (depending on the point $X$ in the spacetime manifold $\mathrm{V}$) \cite{Cattaneo1958, Massa1974a, Massa1974b, Mitskievich1996, JantzenCariniBini1992, A16}. A ``spatial vector" is then defined as a {\it spacetime vector} which is equal to its spatial projection. A full algebra of ``spatial tensors" can be defined in the same way, and also, once a relevant connection has been defined, a spatial tensor analysis \cite{Cattaneo1958, Massa1974a, Massa1974b, Mitskievich1996, JantzenCariniBini1992}. \\

However, it is possible in a general spacetime manifold $\mathrm{V}$ to define a relevant physical space as a {\it 3-D differentiable manifold,} at least {\it locally} in $\mathrm{V}$. To see this, consider a coordinate system or chart: 
\be
\chi : \mathrm{U} \rightarrow \mathbb{R}^4,\ X\mapsto \chi (X)={\bf X}\equiv (x^\mu)\ (\mu =0,...,3),
\ee 
where $\mathrm{U}$ is an open subset of $\mathrm{V}$: the domain of the chart. Then one may define a set of world lines, each of which, $l$, has constant spatial coordinates $a^j$ in the chart $\chi $:
\footnote{\
Note that, if we assume that V is endowed with a Lorentzian metric $\Mat{g}$ whose component $g_{00}$ in the chart $\chi $ verifies $g_{00}>0$ in U, then each among the world lines $l$ is time-like, because in the chart $\chi $ the tangent vector to $l$ has components $\propto (1,0,0,0)$, which may be normed to
\be\label{Vmu}
v^0\equiv \frac{1}{\sqrt{g_{00}}}, \qquad v^j=0.
\ee
}
\be\label{world line}
l_{{\bf a}}=\{X\in \mathrm{U};\ \chi (X)=(x^\mu ) \mathrm{\ is\ such\ that\ }  x^j=a^j\ \mathrm{for}\ j=1,2,3\}.
\ee
Let us suppose for a moment that the chart $\chi $ is in fact a Cartesian coordinate system on the Minkowski spacetime. Then that chart defines an inertial reference frame. In that case, it is clear that, for any event $X$, with  $\chi (X)=(x^\mu )$, the triplet ${\bf x}\equiv (x^j)\ (j=1,2,3)$ defines the spatial position associated in that chart with the event $X$. Note that the data of ${\bf x}$ is equivalent to specifying a unique world line in the ``congruence (\ref{world line})". [By this we mean the set of the world lines (\ref{world line}), when ${\bf a}\equiv (a^j)$ takes any value in $\mathbb{R}^3$ such that the corresponding world line (\ref{world line}) is not empty.] That world line is thus uniquely determined by the event $X$ and may be noted $l(X)$. Events $X'$ that have different values of the time coordinate $x^0$, but that have the same values of the spatial coordinates $x^j \ (j=1,2,3)$, can be said to occur at the same spatial position in the inertial frame as does $X$. Thus, the whole of $l(X)$ is needed. However, each world line in the congruence (\ref{world line}) stays invariant if we change the coordinate system by a purely spatial coordinate change:
\be\label{purely-spatial-change}
x'^0=x^0,\ x'^k=\phi^k((x^j)).
\ee 
It is clear that this transformation leaves us in the same inertial reference frame. With the new chart $\chi '$, the new triplet ${\bf x}'\equiv (x'^k)\equiv \Mat{\phi}({\bf x})$ corresponds to the same spatial position in the inertial frame as does ${\bf x}$ with the first chart $\chi $. And indeed, that world line of the congruence which is defined in the chart $\chi '$ by the data of ${\bf x}'$ is just the same as that world line of the congruence which is defined in the chart $\chi $ by the data of ${\bf x}$. The spatial position of the event $X$ in the inertial frame is therefore most precisely defined by {\it the world line $l(X)$ of the congruence which passes at $X$}. \\

Now note that very little in the foregoing paragraph actually depends on whether or not the chart $\chi $ is a Cartesian coordinate system on the Minkowski spacetime: only the qualification of the reference frame as being inertial depends on that. It is just that we are accustomed to consider a spatial position in an inertial frame in a flat spacetime, and special relativity makes it natural to accept that it is actually the world line $l(X)$ which best represents that spatial position. Hence, consider a general spacetime, and define a congruence of world lines from the data of a coordinate system as in Eq. (\ref{world line}). In the domain of the chart, we may then define the spatial position of an event $X$ as the unique world line $l(X)$ of the congruence (\ref{world line}) which passes at $X$ --- i.e., $l(X)$ is the unique world line of the congruence (\ref{world line}), such that $X\in l(X)$. Thus, {\it the data of a coordinate system on the spacetime defines a three-dimensional space $\mathrm{M}$, of which the \underline{points} (the elements of $\mathrm{M}$) are the \underline{world lines} of the congruence (\ref{world line}) associated with that coordinate system.} 

\subsection{Local reference frame and local space manifold}\label{Local}

The foregoing approach can be used to define precise notions of a reference frame and its unique associated space manifold \cite{A44}. First, the invariance of the congruence (\ref{world line}) under the purely spatial coordinate changes (\ref{purely-spatial-change}) allows one to define a {\it reference frame} as being an equivalence class of charts related by a change (\ref{purely-spatial-change}). More exactly, the following is an equivalence relation between charts which are defined on a given open subspace U of the spacetime manifold V:
\bea\label{R_U}
\chi \mathcal{R}_{\mathrm{U}} \chi '\ \Longleftrightarrow 
[\forall \mathbf{X} \in \chi (\mathrm{U}),\quad \phi^0 (\mathbf{X})=x^0\ \mathrm{and}\  \frac{\partial \phi^k }{\partial x^0}\left(\mathbf{X}\right)=0 \ (k=1,...,3)],
\eea
where $f\equiv \chi '\circ \chi ^{-1}\equiv (\phi^\mu )$ is the transition map, which is defined on $\chi (\mathrm{U})$. Thus a reference frame F is a set of charts defined on the same open domain U and exchanging by a purely spatial coordinate change (\ref{purely-spatial-change}). Then the {\it space manifold} M or $\mathrm{M}_\mathrm{F}$ associated with the reference frame F is defined as the set of the world lines (\ref{world line}). In detail: let $P_S: \mathbb{R}^4 \rightarrow \mathbb{R}^3, {\bf X}\equiv (x^\mu )\mapsto {\bf x}\equiv (x^j )$, be the spatial projection. A world line $l$ is an element of $\mathrm{M}_\mathrm{F}$ iff there is a chart $\chi \in \mathrm{F}$ and a triplet ${\bf x} \in P_S(\chi (\mathrm{U}))$, such that $l$ is the set of {\it all} points $X$ in the domain U, whose spatial coordinates are ${\bf x}$:
\be\label{l-in-M-by-P_S}
l \equiv \{\,X\in \mathrm{U};\ P_S(\chi (X))={\bf x}\,\}.
\ee
(Thus, $l$ is not necessarily a connected set.) It results easily from (\ref{purely-spatial-change}) that (\ref{l-in-M-by-P_S}) holds true then in any chart $\chi' \in \mathrm{F}$, of course with the transformed spatial projection triplet ${\bf x}'= \Mat{\phi }({\bf x})\equiv (\phi ^j({\bf x}))$ \cite{A44}. For a chart $\chi \in \mathrm{F}$, one defines the ``associated chart" as the mapping which associates, with a world line $l\in \mathrm{M}$, the constant triplet of the spatial coordinates of the points $X\in l$:
\be\label{def-chi-tilde}
\widetilde{\chi }: \mathrm{M}\rightarrow \mathbb{R}^3,\quad l\mapsto {\bf x} \mathrm{\ such\ that\ }\forall X \in l,\  P_S(\chi (X))={\bf x}.\\
\ee 
One shows then that the set $\mathcal{T}$ of the subsets $\Omega \subset \mathrm{M}$ such that, 
\be\label{def-Topo}
\forall \chi \in \mathrm{F},\quad \widetilde{\chi }(\Omega)\mathrm{\ is\ an\ open\ set\ in\ }\mathbb{R}^3
\ee
is a topology on $\mathrm{M}$. Finally one shows that the set of the associated charts: $\widetilde{\mathrm{F}}\equiv \{\widetilde{\chi };\ \chi \in \mathrm{F}\}$, is an atlas on the topological space $(\mathrm{M},\mathcal{T})$, hence defines a structure of {\it differentiable manifold} on $\mathrm{M}$ \cite{A44}.  {\it Thus the space manifold $\mathrm{M}_\mathrm{F}$ is browsed by precisely the triplet ${\bf x}\equiv (x^j)$ made with the spatial projection of the spacetime coordinates ${\bf X}\equiv (x^\mu )\equiv \chi (X)$ of a chart $\chi \in \mathrm{F}$, see Eq. (\ref{def-chi-tilde}).}\\

Using these results, one may define the space of quantum-mechanical states, for a given reference frame F in a given spacetime $(\mathrm{V},\Mat{g})$, as being the set $\mathcal{H}$ of the square-integrable functions defined on the corresponding space manifold M \cite{A43}. One may also define the full algebra of spatial tensors: the pointwise algebra is defined simply as the tensor algebra of the tangent vector space TM$_x$ to the space manifold M at some arbitrary point $x \in \mathrm{M}$ \cite{A47}.

\subsection{Goal and summary}\label{Summary}

Thus, by defining a reference frame as a set F of charts that all have the same open domain U and that exchange by a purely spatial coordinate change (\ref{purely-spatial-change}), one can then define the associated space manifold $\mathrm{M}_\mathrm{F}$ as the set of the world lines (\ref{l-in-M-by-P_S}) \cite{A44}. These definitions are relevant to physical applications \cite{A43,A47}. However, they apply to a {\it parametrizable} domain U of the spacetime manifold V, i.e., to an open set U, such that at least one regular chart can be defined over the whole of U. Since the manifold V itself as a whole is in general not parametrizable, a reference frame is in general only a local one, and so the associated space manifold does not look ``maximal". {\it The aim of the present work} is to define {\it global} reference fluids, to associate with any of them a {\it global physical space,} and to {\it link} these concepts with the formerly defined {\it local concepts.} As the ``spacetime", we consider a differentiable manifold $\mathrm{V}$ having dimension $N+1$, thus $N$ is the dimension of the ``space" manifold to be defined. We define a global reference fluid by the data of a {\it non-vanishing} global vector field $v$ on $\mathrm{V}$. We do not need that $N=3$, nor that $\mathrm{V}$ be endowed with a Lorentzian metric $\Mat{g}$ for which $v$ be a time-like vector field. This was already true for the former ``local" work \cite{A44}. Note, however, that a time-like vector field on a Lorentzian manifold $(\mathrm{V},\Mat{g})$ is non-vanishing; and that, if a Lorentzian manifold $(\mathrm{V},\Mat{g})$ is time-oriented, which indeed is usually required for a spacetime, then by definition there exists at least one global time-like vector field on $\mathrm{V}$. We define the ``global space" associated with $v$ as the set $\mathrm{N}_v$ of the maximal integral curves (or ``orbits") of $v$. To reach our goal, we take the following steps: \\

\noi ({\it a}) Section \ref{AdaptedCharts} studies when a given vector field $v$ on a differentiable manifold V is such that there locally exists charts of V which are {\it adapted} \cite{Cattaneo1958} to the congruence associated with $v$; i.e., charts in which the ``spatial" position ${\bf x}\equiv (x^j) \ (j=1,...,N)$ is constant on any orbit $l$ of $v$; see \hyperref[Definition 1]{Definition 1}. We need also that the mapping $l\mapsto {\bf x}$ be injective. The desired situation is defined by \hyperref[Proposition 2]{Proposition 2}. According to a \hyperref[Conjecture]{transversality argument}, this situation should be attainable, in general, if $v$ does not vanish and each of its orbits is non-periodic and is closed in $\mathrm{V}$. In Subsect. \ref{Transverse}, two theorems of transversality and another theorem pertaining to differential topology allow us to formalize that argument in \hyperref[Theorem 4]{Theorem 4}. This justifies us in introducing a notion of ``normal" vector field by \hyperref[Normal]{Definition 2}, which ensures the local existence of $v$-adapted charts through \hyperref[Theorem 5]{Theorem 5}.\\

\noi ({\it b}) For any $v$-adapted chart $\chi $, one may define the mapping $\bar{\chi }$ which associates with any orbit $l$ the constant spatial position ${\bf x}$. We show in Section \ref{N_v manifold} that, using the set $\mathcal{A}$ of the injective mappings $\bar{\chi }$, one can endow the global orbit set $\mathrm{N}_v$ with a topology $\mathcal{T}'$ for which this set is an atlas of charts. See \hyperref[Theorem 6]{Theorem 6}. Thus, when $\mathcal{T}'$ is metrizable and separable, we do have a canonical structure of differentiable manifold on the orbit set $\mathrm{N}_v$, for which the mappings $\bar{\chi }$ defined by Eq. (\ref{Def chi bar}) are local charts on the ``space" manifold $\mathrm{N}_v$. I.e., if $N=3$, {\it also the global space $\mathrm{N}_v$ is browsed (locally) by precisely the triplet ${\bf x}\equiv (x^j)$ made with the spatial projection of the spacetime coordinates ${\bf X}\equiv (x^\mu )\equiv \chi (X)$ of a $v$-adapted chart $\chi $.}\\

\noi ({\it c}) In Section \ref{Local_vs_global} we establish the link with the previously defined space manifold $\mathrm{M}_\mathrm{F}$, associated with a given local reference frame F --- defined as an equivalence class of charts for the relation (\ref{R_U}). We show in \hyperref[Theorem 7]{Theorem 7} that, when the charts belonging to $\mathrm{F}$ are $v$-adapted and with the mapping $l\mapsto {\bf x}$ being injective, then $\mathrm{M}_\mathrm{F}$ is naturally identified with an open subset of $\mathrm{N}_v$. \\

\noi The definitions of a local reference frame F and the corresponding ``space" manifold $\mathrm{M}_\mathrm{F}$ do not need any metrical structure on the ``spacetime" manifold $\mathrm{V}$  \cite{A44}. Just the same can be said for the definition of the global ``space" manifold $\mathrm{N}_v$, beyond the very fact that $\mathrm{V}$ should have a metrizable topology. 

\section{Local existence of adapted charts}\label{AdaptedCharts}

\subsection{Definitions}\label{Definitions}

Let $\mathrm{V}$ be an $(N+1)-$dimensional real differentiable manifold, with $N\geq 1$.
\footnote{\label{Def Dif Manifold}\,
We understand ``differentiable manifold" as a topological space endowed with an atlas of compatible charts, hence with the corresponding equivalence class of compatible atlases --- with the restriction that that space should be metrizable and separable \cite{DieudonneTome3}. 
}
We consider a given global, smooth vector field $v$ on $\mathrm{V}$. A continuously differentiable ($\mathcal{C}^1$) mapping $C$ from an open interval I of $\mathbb{R}$ into V defines an integral curve of $v$ iff $\frac{\dd C}{\dd s}=v(C(s))$ for $s\in \mathrm{I}$. For any $X\in \mathrm{V}$, let $C_X$ be the solution of 
\be\label{Def C}
\frac{\dd C}{\dd s}=v(C(s)),\qquad C(0)=X
\ee
for which the open interval I is {\it maximal,} and denote this maximal interval by $\mathrm{I}_X$ \cite{DieudonneTome4}. That is, $\mathrm{I}_X$ is an open interval defined as the union of all open intervals I, each containing 0, in which a solution of (\ref{Def C}) is defined. The solution $C_X$ is defined on $\mathrm{I}_X$ and is unique \cite{DieudonneTome4}. Let $s \in \mathrm{I}_X$ and set $Y=C_X(s)$. These definitions imply easily \cite{DieudonneTome4} that 
\be\label{C_X vs C_Y}
\mathrm{I}_Y=\mathrm{I}_X-s\quad \mathrm{and\quad} \forall t\in \mathrm{I}_Y,\ C_{Y}(t)=C_X(s+t).
\ee
For any $X\in \mathrm{V}$, call the {\it range} $\,l_X\equiv C_X(\mathrm{I}_X)\subset \mathrm{V}$ the ``maximal integral curve at $X$". From (\ref{C_X vs C_Y}), ``$l_X$ does not depend on the point $X\in l_X$": if $Y\in l_X$, then $l_Y=l_X$. We define the set of the maximal integral curves (or orbits) of $v$:
\be\label{Def N_v}
\mathrm{N}_v\equiv \{l_X;\ X\in \mathrm{V}\}.
\ee
Once endowed with further structure, $\mathrm{N}_v$ will be the global space manifold associated with the global vector field $v$ (when the latter is non-vanishing and obeys another assumption). Note that, if the set $\mathrm{U}\subset \mathrm{V}$ is not empty, then the following subset of $\mathrm{N}_v$: 
\be\label{Def D_U}
\mathrm{D}_\mathrm{U}\equiv \{l\in \mathrm{N}_v;\ l\cap \mathrm{U} \ne \emptyset \}
\ee
is not empty. Indeed, for any $X\in \mathrm{U}$, the world line $l\equiv l_X$ belongs to $\mathrm{D}_\mathrm{U}$. Let $ P_S: \mathbb{R}^{N+1} \rightarrow \mathbb{R}^N,\quad (x^\mu )\mapsto (x^j)\ (\mu =0,...,N;\ j=1,...,N)$ be the ``spatial" projection. 

\paragraph{Definition 1.}\label{Definition 1} {\it A mapping $\chi:\mathrm{U}\rightarrow\mathbb{R}^{N+1}$ with $\mathrm{U}\subset \mathrm{V}$ is said ``$v$-adapted" iff for any $l\in \mathrm{D}_\mathrm{U}$, there exists ${\bf x} \in \mathbb{R}^N$ such that}
\be\label{Def adapted}
\forall Y\in l\cap \mathrm{U},\quad P_S(\chi (Y))={\bf x}.
\ee
If Eq. (\ref{Def adapted}) is verified by some world line $l\in \mathrm{N}_v$, then necessarily $l\in \mathrm{D}_\mathrm{U}$, and ${\bf x}$ is obviously unique. Thus, for any $v$-adapted mapping $\chi$, the mapping 
\be\label{Def chi bar}
\bar{\chi}: \mathrm{D}_\mathrm{U} \rightarrow \mathbb{R}^N,\quad l\mapsto {\bf x} \mathrm{\ \,such\ that\ (\ref{Def adapted})\ is\ verified}
\ee
is well defined. In Section \ref{N_v manifold}, we will endow the set $\mathrm{N}_v$ with first a topology and then a structure of differentiable manifold, for which the charts (of $\mathrm{N}_v$) will be mappings $\bar{\chi}$, where $\chi $ is a $v$-adapted {\it chart} of $\mathrm{V}$. Since any chart is in particular a one-to-one mapping, we need to restrict ourselves to $v$-adapted charts $\chi $ such that the associated  mapping $\bar{\chi}$ is injective. Thus, we define that a $v$-adapted chart $\chi $ is {\it ``nice"} iff $\bar{\chi}$ is injective on $\mathrm{D}_\mathrm{U}$, with $\mathrm{U}\subset \mathrm{V}$ the (open) domain of $\chi $. 

\subsection{Straightening-out vs $v$-adapted charts}

In the remainder of this section, we investigate whether there exist nice $v$-adapted charts in the neighborhood of any point $X_0\in \mathrm{V}$. If the vector field $v$ does not vanish, a well-known theorem applies at any point $X_0\in \mathrm{V}$:

\paragraph{Straightening-out theorem}\label{Straightening-out theorem} (e.g. \cite{AbrahamMarsdenRatiu}). {\it Let $v$ be a vector field of class $\mathcal{C}^\infty$ defined on $\mathrm{V}$. Suppose that at $X_0\in \mathrm{V}$ we have $v(X_0) \ne 0$. There is a ``straightening-out chart" $\chi $ defined on an open neighborhood $\mathrm{U}$ of $X_0$, i.e. $\chi $ is such that:}

\bi
\item (i) $\chi (\mathrm{U}) =\mathrm{I}\times \Omega$, $\ \mathrm{I}=\,]-a,+a[\,$,$\ a\ne 0,\quad \Omega\ \mathrm{open\ set\ in\ }\mathbb{R}^N$. 

\item (ii) {\it For any ${\bf x}\in \Omega $, $\chi ^{-1}(\mathrm{I}\times \{{\bf x}\})$ is an integral curve of $v$.}

\item (iii) {\it In $\mathrm{U}$, we have $v=\partial _0$, where $(\partial _\mu )$ is the natural basis associated with the chart $\chi $.}

\ei 

\noi However, the direct link with the notion of a $v$-adapted chart works in the wrong direction:

\paragraph{Proposition 0.}\label{Proposition 0} {\it Let $(\chi, \mathrm{U})$ be a $v$-adapted chart.} (i) {\it We have $v_{\mid \mathrm{U}}=f\partial _0$ with $f: \mathrm{U}\rightarrow \mathbb{R}$ a smooth function.} (ii) {\it Given any point $X\in \mathrm{U}$ such that $v(X)\ne 0$, one may obtain a straightening-out chart $(\chi', \mathrm{U}')$, with $\mathrm{U}'\subset \mathrm{U}$ being an open neighborhood of $X$, by changing merely the $y^0$ coordinate.} \\

\noi {\it Proof.} (i) To say that $(\chi, \mathrm{U})$ is $v$-adapted means, according to \hyperref[Definition 1]{Definition 1}, that for any given $X\in \mathrm{U}$, we have for any $Y\in  l_X\cap \mathrm{U}$: $P_S(\chi (Y))=P_S(\chi (X))$. In the coordinates $\chi (X)=(x^\mu )\equiv {\bf X},\ \chi (Y)=(y^\mu )\equiv {\bf Y}\ (\mu =0,...,N)$, the latter rewrites as 
\be\label{y^j=x^j}
y^j=x^j\quad (j=1,...,N).
\ee
On the other hand, remembering the definition of $l_X$ as a (maximal) integral curve of $v$, Eq. (\ref{Def C}) and below, let $\mathrm{J}$ be the connected component of $0$ in $\mathrm{I}'_{X \mathrm{U}}\equiv \{s\in \mathrm{I}_X;\, C_X(s)\in \mathrm{U}\}$: $\mathrm{J}$ is an open interval containing $0$ and we have $l'\equiv C_X(\mathrm{J})\subset l_X\cap \mathrm{U}$. Let us denote this as $l'=\{Y(s);s\in \mathrm{J}\}$, with
\be\label{l' in chi}
\frac{\dd y^\mu }{\dd s}=v^\mu({\bf Y}(s))\quad  (\mu =0,...,N), \quad {\bf Y}(0)={\bf X},
\ee
$v^\mu =v^\mu({\bf Y})$ being the components of $v$ in the chart $\chi $. It follows from (\ref{y^j=x^j}) and (\ref{l' in chi}) that for $s\in \mathrm{J}$ we have $v^j({\bf Y}(s))=0\ (j=1,...,N)$, i.e. $v(Y(s))=v^0({\bf Y}(s))\partial _0(Y(s))$, thus in particular $v(X)=v^0({\bf X})\partial _0(X)\equiv f(X)\partial _0(X)$. Since this is true at any point $X\in \mathrm{U}$, our statement (i) is proved. \\

(ii) If we leave the coordinates $y^j\ (j=1,...,N)$ unchanged: $\chi' (Y)\equiv {\bf Y}'=(g({\bf Y}),(y^j))$ where ${\bf Y}=(y^0,(y^j))=\chi (Y)\quad (Y\in \mathrm{U}'\subset \mathrm{U})$, then the components of $v$ in the new chart $(\chi', \mathrm{U}')$ are $v'^j=v^j=0\ (j=1,...,N)$ and $v'^0({\bf Y}')=\frac{\partial g}{\partial y^0}v^0({\bf Y})=\frac{\partial g}{\partial y^0}f(Y)$. The latter must be equal to $1$ for a straightening-out chart. Since $v_{\mid \mathrm{U}}=f\partial _0$ and $v(X)\ne 0$, we have $f(Y)\ne 0$ when $Y$ is  in some neighborhood $\mathrm{U}''\subset \mathrm{U}$ of $X$. Hence, we get $v'^0=1$ in $\mathrm{U}'$ if we take $\mathrm{U}'=\chi ^{-1}(\mathrm{B})$ with $\mathrm{B}=]x^0-r,x^0+r[\times ...\times ]x^N-r,x^N+r[\subset \chi (\mathrm{U}'')$ and define $y'^0 \equiv g({\bf Y})=\int _{x^0} ^{y^0}\dd u/f(u,(y^j))$ for ${\bf Y}\in \mathrm{B}$. Property (i) in the straightening-out theorem is then got by a mere shift $y'^0\hookrightarrow y'^0+\delta $, and its Property (ii) is a straightforward consequence of Properties (i) and (iii). \hfill $\square$\\

Conversely, if $(\chi, \mathrm{U}) $ is a straightening-out chart, i.e. it fulfills conditions (i) to (iii) of the \hyperref[Straightening-out theorem]{theorem above}, then let ${\bf x}\in \Omega $ and set $l'\equiv \chi ^{-1}(\mathrm{I}\times \{{\bf x}\})$. From (ii), $l'$ is an integral curve of $v$ and, since $l'\subset \mathrm{U}$ by construction, we have (\ref{Def adapted}):
\be\label{Def adapted'}
\forall Y\in l'\cap \mathrm{U},\quad P_S(\chi (Y))={\bf x}.
\ee
Hence, at first sight, it might seem that $\chi $ is a $v$-adapted chart. However, let $X=\chi ^{-1}(s,{\bf x})\in l'$ and let $l\equiv l_X \in \mathrm{N}_v$ be the {\it maximal} integral curve of $v$ passing at $X$. We have $l'\subset (l\cap \mathrm{U})$ since $l'$ is an integral curve of $v$ that is included in $\mathrm{U}$ and that passes at $X$. But nothing guarantees that $l'\supset  (l\cap \mathrm{U})$: the intersection $l\cap \mathrm{U}$ may contain other arcs, say $l'_1\equiv \chi ^{-1}(\mathrm{I}\times \{{\bf x}_1\})$ with ${\bf x}_1 \in \Omega $ and ${\bf x}_1 \ne {\bf x}$. In such a case, the straightening-out chart $\chi $ is not $v$-adapted, since for $Y\in l\cap \mathrm{U}$, $P_S(\chi (Y))$ may take different values ${\bf x}$, ${\bf x}_1$, ... \\

As we already noted, Property (ii) in the straightening-out theorem follows easily from Properties (i) and (iii). More generally, if a chart $(\chi ,\mathrm{U})$ satisfies Property (iii), i.e.  $v=\partial _0$ in  $\mathrm{U}$, and if $\chi (\mathrm{U})$ contains a set $\mathrm{I}\times\{{\bf x}\}$, with ${\bf x}\in \mathbb{R}^N$ and $\mathrm{I}$ an open interval, then Property (ii) applies for this ${\bf x}\in \mathbb{R}^N$. Later on, we will need that the boundary of the open set $\mathrm{U}$ be a smooth hypersurface, hence we must consider charts for which (iii) is true, but not (i).

\paragraph{Proposition 1.}\label{Proposition 1} {\it Let $(\chi, \mathrm{U})$ be a chart such that $v=\partial _0$ in  $\mathrm{U}$. Assume there is an open subset $\Omega \subset \mathbb{R}^N$ such that
\be\label{chi(U)}
\chi (\mathrm{U})=\bigcup _{{\bf x}\in \Omega }\mathrm{I}_{\bf x} \times \{{\bf x}\},
\ee
where $\mathrm{I}_{\bf x}$ (${\bf x}\in \Omega $) are open intervals.} \\
\noi (i) {\it In order that the chart $(\chi, \mathrm{U})$ be $\,v\,$-adapted, it is necessary and sufficient that}
\be\label{Straightening->v-adapted}
\forall X \in \mathrm{U}, \ \chi(l_X\cap \mathrm{U}) \mathrm{\ has\ the\ form\ }\mathrm{I}_{\bf x} \times \{{\bf x}\} \mathrm{\ for\  some\ }{\bf x}\in \Omega .
\ee
(ii) {\it Moreover, if that is the case, then the $\,v\,$-adapted chart $(\chi, \mathrm{U}) $ is nice.}\\

\noi {\it Proof.} (i) If Condition (\ref{Straightening->v-adapted}) is satisfied, let $l\in \mathrm{D}_\mathrm{U}$. Thus, $l\cap \mathrm{U} \ne \emptyset$, so let $X\in l\cap \mathrm{U}$, hence $l=l_X$. Let ${\bf x}\in \Omega $ be given by (\ref{Straightening->v-adapted}) for precisely the maximal integral curve $l=l_X$. For any $Y\in l\cap \mathrm{U}$, we have thus $\chi (Y)=(s,{\bf x})$ for some $s\in \mathrm{I}_{\bf x}$. Hence, we have (\ref{Def adapted}). Therefore, according to \hyperref[Definition 1]{Definition 1}, $\chi$ is $\,v\,$-adapted. Conversely, assume that  $(\chi, \mathrm{U})$ is $\,v\,$-adapted. Let $X\in \mathrm{U}$ and set $l\equiv l_X$. Further, let ${\bf x}$ be given by (\ref{Def adapted}). That is, any point ${\bf Y}\in \chi (l\cap \mathrm{U})$ has the form $(s,{\bf x})$ for some $s\in \mathbb{R}$. Since moreover $\chi(\mathrm{U})$ has the form (\ref{chi(U)}), we have also ${\bf Y}=(t,{\bf x}')$ for some ${\bf x}'\in \Omega $ and some $t\in \mathrm{I}_{{\bf x}'}$, so ${\bf x}={\bf x}'$ and $s=t\in \mathrm{I}_{\bf x}$. Hence $\chi(l\cap \mathrm{U}) \subset \mathrm{I}_{\bf x} \times \{{\bf x}\}$. But also $\mathrm{I}_{\bf x}\times \{{\bf x}\}\subset \chi (l\cap \mathrm{U})$. Indeed, setting $l'\equiv \chi ^{-1}(\mathrm{I}_{\bf x}\times \{{\bf x}\})$, we have $l'\subset (l\cap \mathrm{U})$ because, as noted before the statement of this Proposition 1, $l'$ is an integral curve of $v$ that is included in $\mathrm{U}$ and that passes at $X$. \\

\noi (ii) Assuming that the chart $(\chi, \mathrm{U})$ obeys Condition (\ref{Straightening->v-adapted}) [and hence, by (i), is a $\,v\,$-adapted chart], let us show that the mapping $\bar{\chi}$ defined by (\ref{Def chi bar}) is injective. Thus, let $l,l'\in \mathrm{N}_v$, assume that both intersect $\mathrm{U}$, and let ${\bf x}$ and ${\bf x}'$ be the images of $l$ and $l'$ by $\bar{\chi}$. This means, according to the definition (\ref{Def chi bar}), that we have (\ref{Def adapted}), and similarly
\be\label{Def adapted''}
\forall Y\in l'\cap \mathrm{U},\quad P_S(\chi (Y))={\bf x}'.
\ee
From Condition (\ref{Straightening->v-adapted}), we get thus:
\be
\chi(l\cap \mathrm{U})=\mathrm{I}_{{\bf x}}\times\{{\bf x}\}\quad \mathrm{and}\quad \chi(l'\cap \mathrm{U})=\mathrm{I}_{{\bf x}'}\times\{{\bf x}'\}.
\ee
Therefore, if ${\bf x}={\bf x}'$, it is clear that $l=l'$. The proof is complete.  \hfill $\square$ \\

\noi The condition that the open set $\mathrm{A}\equiv \chi (\mathrm{U)}\subset \mathbb{R}^{N+1}$ have the form (\ref{chi(U)}) is fulfilled, in particular, if $\mathrm{A}$ is convex. So it is fulfilled if one restricts the chart $\chi $ to $\chi ^{-1}(\mathrm{A})$ with $\mathrm{A}$ a convex open subset of $\chi (\mathrm{U}_0)$. Unfortunately, what we have in a rather general situation is the following result, which does not ensure the applicability of Proposition 1:

\paragraph{Theorem 0}\label{Theorem 0} \cite{DieudonneTome4-2}. {\it Suppose that some maximal integral curve $l$ of the $\mathcal{C}^\infty$ vector field $v$ is closed in $\mathrm{V}$ and is not reduced to a point. Since $v$ then does not vanish on $l$, let $\chi: \mathrm{U} \rightarrow \mathrm{I}\times \Omega$ be a straightening-out chart in an open neighborhood $\mathrm{U}$ of some point of $l$. Then we have $\chi (l\cap \mathrm{U})=\mathrm{I}\times \mathrm{E}$, where $\mathrm{E}$ is a closed countable subset of $\Omega $, and any point $X\in l\cap \mathrm{U}$ is ``transversally isolated", i.e. ${\bf x}\equiv P_S(\chi (X))$ is isolated in $\mathrm{E}$.}\\

\noi (Note that the set $\mathrm{E}$ is closed in $\Omega $, which is an open set in $\mathbb{R}^N$. So $\mathrm{E}$ is not necessarily closed in $\mathbb{R}^N$.) Assume that some straightening-out chart $(\chi,\mathrm{U}) $ is such that, for any $X\in \mathrm{U}$, the maximal integral curve $l_X$ is closed in $\mathrm{V}$. From Point (i) in \hyperref[Proposition 1]{Proposition 1}, we get that this is a $v$-adapted chart iff, in addition, for any point $X \in \mathrm{U}$, the countable closed subset $\mathrm{E}_X$ of $ \Omega $, whose existence is ensured by Theorem 0 for the curve $l_X$, is actually reduced to a point. We can further characterize this desired situation:

\paragraph{Proposition 2.}\label{Proposition 2} {\it Let $\chi: \mathrm{U}_0 \rightarrow \mathrm{I}\times \Omega_0$ be a straightening-out chart for the $\mathcal{C}^\infty$ vector field $v$. Let $\mathrm{U}\subset \mathrm{U}_0$ be an open set such that $\chi (\mathrm{U})$ has the form (\ref{chi(U)}). Assume that, for any $X\in \mathrm{U}$, the maximal integral curve $l_X$ is closed in $\mathrm{V}$.} \\
\noi (i) {\it For any given $X \in \mathrm{U}$, set  $\chi (X)=(s,{\bf x})$. The connected component $\lambda ''$ of $(s,{\bf x})$ in $\lambda \equiv \chi (l_X\cap \mathrm{U})$ is equal to $\lambda' \equiv \mathrm{I}_{\bf x}\times\{{\bf x}\}$.} \\
\noi (ii) {\it In order that $(\chi,\mathrm{U}) $ be a $v$-adapted chart, it is necessary and sufficient that, for any $X\in \mathrm{U}$, the intersection $l_X\cap \mathrm{U}$ be a connected set.}\\
\noi (iii) {\it Let $\mathrm{W}$ be an open subset of $\mathrm{U}$ such that, for any $X\in \mathrm{W}$, $l_X\cap \mathrm{U}$ be a connected set. Then the restriction $(\chi,\mathrm{W}) $ is a nice $v$-adapted chart.}\\

\noi {\it Proof.} (i) Since $\mathrm{I}_{\bf x}$ is an interval, the set $\lambda '\equiv \mathrm{I}_{\bf x}\times\{{\bf x}\}$  is connected. To show that $\lambda ''=\lambda '$ is always true, we note first that, $\chi (\mathrm{U})$ having the form (\ref{chi(U)}), $\chi (X)=(s,{\bf x})$ is such that $s\in \mathrm{I}_{\bf x}$. Hence $(s,{\bf x})\in \lambda '$, and since $\lambda '$ is connected, we have $\lambda '\subset \lambda ''$. To prove that in fact $\lambda '= \lambda ''$, we will show that $\lambda '$ is open and closed in $\lambda ''$. We know from \hyperref[Theorem 0]{Theorem 0} that $\lambda_0\equiv \chi (l_X\cap \mathrm{U}_0)$ has the form $\lambda_0 = \mathrm{I}\times \mathrm{E}$, where $\mathrm{E}$ is a subset of $\Omega_0 $ having only isolated points. Thus, ${\bf x}$ being isolated in  $\mathrm{E}$, let $r>0$ be such that $\mathrm{B}\cap \mathrm{E}=\{{\bf x}\}$, where $\mathrm{B}\equiv \mathrm{B}({\bf x},r)$ is the open ball of radius $r$ in $\mathbb{R}^N$, centered at ${\bf x}$. Hence, we have 
\be
(\mathrm{I}_{\bf x}\times \mathrm{B})\cap \lambda ''\subset (\mathrm{I}_{\bf x}\times\mathrm{B} )\cap \lambda _0=(\mathrm{I}_{\bf x}\times\mathrm{B} )\cap(\mathrm{I}\times \mathrm{E})=(\mathrm{I}_{\bf x}\cap\mathrm{I})\times (\mathrm{B} \cap\mathrm{E})=\mathrm{I}_{\bf x}\times\{{\bf x}\}\equiv \lambda '.
\ee
We have also $\lambda '\subset (\mathrm{I}_{\bf x}\times \mathrm{B})\cap \lambda ''$, because $\lambda '\subset \mathrm{I}_{\bf x}\times \mathrm{B}$ and $\lambda '\subset \lambda ''$. So $\lambda '=(\mathrm{I}_{\bf x}\times \mathrm{B})\cap \lambda ''$ is an open subset of $\lambda ''$. On the other hand, if a sequence of points of $\lambda '$ tends towards a limit in $\lambda ''$, say $(s_n,{\bf x})\rightarrow (s',{\bf y})\in \lambda ''$, then ${\bf y}={\bf x}$ and, since $\lambda ''\subset \chi (\mathrm{U})$ which is given by (\ref{chi(U)}), we have $(s',{\bf x})\in (\mathbb{R}\times \{{\bf x}\})\cap  \chi (\mathrm{U})=\mathrm{I}_{\bf x}\times\{{\bf x}\}\equiv \lambda '$, hence the limit $(s',{\bf x})$ is in $\lambda '$, so that $\lambda '$ is a closed subset of $\lambda ''$. Being non-empty and an open-and-closed subset of the connected set $\lambda ''$, $\lambda '$ is equal to $\lambda ''$. \\

\noi (ii) Since $\chi : \mathrm{U}_0 \rightarrow \mathrm{I}\times \Omega_0$ is a bicontinuous mapping, it is of course equivalent to say that $l_X\cap \mathrm{U}$ or $\lambda _X\equiv \chi (l_X\cap \mathrm{U})$ is connected. Therefore, (ii) follows immediately from (i) and from Statement (i) in \hyperref[Proposition 1]{Proposition 1}. \\

\noi (iii) Let $X$ be any point in $ \mathrm{W}$ and set $\chi (X)=(s,{\bf x})$. Since $l_X\cap \mathrm{U}$ is connected, we have $\chi (l_X\cap \mathrm{U})=\mathrm{I}_{\bf x}\times\{{\bf x}\} $ from (i). Thus, $P_S(\chi (Y))={\bf x}$ is true for any $Y\in l_X\cap \mathrm{U}$, hence a fortiori for any $Y\in l_X\cap \mathrm{W}$. Hence the chart $(\chi, \mathrm{W})$ is $v$-adapted. In a similar way, it is easy to adapt the proof of Statement (ii) in \hyperref[Proposition 1]{Proposition 1} to conclude that the $v$-adapted chart $(\chi, \mathrm{W})$ is nice. \hfill $\square$

\paragraph{Proposition 3.}\label{Proposition 3} {\it Let $\chi: \mathrm{U}_0 \rightarrow \mathrm{I}\times \Omega_0$ be a straightening-out chart for the $\mathcal{C}^\infty$ vector field $v$ and assume that, for some $X\in \mathrm{U}_0$, the maximal integral curve $l_X$ is closed. Let $\mathrm{E}$ be the closed countable subset of $\Omega_0 $ given by \hyperref[Theorem 0]{Theorem 0} for the maximal curve $l\equiv l_X$, thus $\chi (l_X\cap \mathrm{U}_0) = \mathrm{I}\times \mathrm{E}$. Set ${\bf x}\equiv P_S(\chi (X))$ and, as ensured by \hyperref[Theorem 0]{Theorem 0}, let $\Omega \subset \Omega_0 $ be any open neighborhood of ${\bf x}$ such that $\mathrm{E}\cap \Omega =\{{\bf x}\}$. Let $\mathrm{U}\subset \mathrm{U}_0$ be any open set such that $\chi (\mathrm{U})$ has the form (\ref{chi(U)}) with this set $\Omega $. Then  $l_X\cap \mathrm{U}$ is connected, $\chi (l_X\cap \mathrm{U})=\mathrm{I}_{\bf x}\times\{{\bf x}\}$.}\\

\noi {\it Proof.} Since $\chi: \mathrm{U}_0 \rightarrow \mathrm{I}\times \Omega_0$ is a bijection, we have
\bea
\chi (l_X\cap \mathrm{U}) & = &\chi ((l_X\cap \mathrm{U}_0)\cap \mathrm{U})= \chi (l_X\cap \mathrm{U}_0)\cap \chi (\mathrm{U})\\
 \nonumber & = &(\mathrm{I}\times\mathrm{E})\cap \left(\bigcup _{{\bf x}'\in \Omega }\mathrm{I}_{{\bf x}'} \times \{{\bf x}'\}\right)=\bigcup _{{\bf x}'\in \Omega }(\mathrm{I}\cap \mathrm{I}_{{\bf x}'}) \times (\mathrm{E}\cap  \{{\bf x}'\})\\
\nonumber & =& \mathrm{I}_{\bf x}\times\{{\bf x}\}.  \hspace{80mm}  \square
\eea
\vspace{3mm}

\subsection{Intersections of straight lines with inverse images under the flow}\label{Straight inter}

Recall that the flow of the global vector field $v$ on $\mathrm{V}$ is the mapping $F: \mathcal{D}\rightarrow \mathrm{V}, \quad (s,X)\mapsto F(s,X)\equiv C_X(s)$, where $\mathcal{D}$ is the domain of the flow $F$:
\be\label{Def D}
\mathcal{D} \equiv \bigcup _{X\in \mathrm{V}} \mathrm{I}_X\times \{X\} \subset \mathbb{R}\times \mathrm{V}.
\ee
If $v$ is $\mathcal{C}^q$ \ ($q \ge 1$, possibly $q =\infty )$, $\mathcal{D}$ is an open set in $\mathbb{R}\times \mathrm{V}$, moreover $F$ is $\mathcal{C}^q$ on $\mathcal{D}$ \cite{DieudonneTome4, AbrahamMarsdenRatiu}. 

\paragraph{Proposition 4.}\label{Proposition 4} {\it Let $\mathrm{U}$ be an open subset of $\mathrm{V}$.} (i) {\it For any $X\in \mathrm{V}$, we have 
\be
l_X\cap \mathrm{U}=F_X(\mathrm{I}'_{X \mathrm{U}})=F(\mathrm{I}'_{X \mathrm{U}}\times \{X\}), 
\ee
where $F_X\equiv F(.,X)=C_X$ is defined on $\mathrm{I}_X \subset \mathbb{R}$, and where}
\be\label{Def I'_X}
\mathrm{I}'_{X \mathrm{U}}\equiv F_X^{-1}(\mathrm{U})=\{s\in \mathrm{I}_X;\, F(s,X)\in \mathrm{U}\}.
\ee
(ii) {\it Further, we have for any $X\in \mathrm{V}$:
\be\label{I'_X = R inter D_U}
\mathrm{I}'_{X \mathrm{U}}\times \{X\}=(\mathbb{R}\times \{X\})\cap \mathcal{D}_\mathrm{U}=(\mathrm{I}_X\times \{X\})\cap \mathcal{D}_\mathrm{U},
\ee
where}
\be\label{Dcal_U}
\mathcal{D}_\mathrm{U}\equiv F^{-1}(\mathrm{U})=\bigcup _{X\in \mathrm{V}} \mathrm{I}'_{X \mathrm{U}}\times \{X\}.
\ee
(iii) {\it  Assume that $v$ is $\mathcal{C}^\infty $ and that, for some $X\in \mathrm{V}$, $l_X$ is closed in $\mathrm{V}$ and not reduced to a point. Then ({\it a}) $l_X$ is a submanifold of $\mathrm{V}$ and the mapping $F_X: \mathrm{I}_X\rightarrow l_X$ is a local diffeomorphism at any point $s\in \mathrm{I}_X$. ({\it b}) Assume moreover that $F_X$ is non-periodic. Then it is a (global) diffeomorphism of $\mathrm{I}_X$ onto $l_X$. The connected components of $l_X\cap \mathrm{U}$ are the images by $F_X$ of the connected components of $\mathrm{I}'_{X \mathrm{U}}\equiv F_X^{-1}(\mathrm{U})=F_X^{-1}(l_X\cap \mathrm{U})$, which are open intervals of $\mathbb{R}$. In particular, in order that $l_X\cap \mathrm{U}$ be connected, it is necessary and sufficient that $\mathrm{I}'_{X \mathrm{U}} $ be an open interval.}\\

\noi {\it Proof.} Points (i) and (ii) follow immediately from the definitions. Let us prove Point (iii). ({\it a}) Since the maximal integral curve $l_X$ is closed in $\mathrm{V}$, this is a submanifold of $\mathrm{V}$ \cite{DieudonneTome4-2}. Since $l_X$ is not reduced to a point, the vector field $v$ does not vanish on $l_X$. [If $v(Y)=0$ for some $Y\in l_X$, then we have $\mathrm{I}_Y=\mathbb{R}$ and $C_{Y}(s)=Y\ \forall s\in \mathbb{R}$ from the uniqueness of the maximal integral curve, hence $X=Y$ from the translation invariance (\ref{C_X vs C_Y}).] Therefore, $\frac{\dd F_X}{\dd s}=v(F_X(s))\ne 0$ implies that $F_X$ is a local diffeomorphism, at any point $s\in \mathrm{I}_X$, between the one-dimensional manifolds $\mathrm{I}_X$ and $l_X$. ({\it b}) Since $l_X$ is closed in $\mathrm{V}$, not reduced to a point, and since $F_X$ is non-periodic, it follows that $F_X$ is injective \cite{DieudonneTome4-2}. In view of ({\it a}) and since $F_X$ is surjective by the definition of $l_X$, it is a diffeomorphism of $\mathrm{I}_X$ onto $l_X$. Hence, the connected components $l_j$ of $l_X\cap \mathrm{U}$ are the images by $F_X$ of the connected components of the open set $\mathrm{I}'_{X \mathrm{U}}\equiv F_X^{-1}(\mathrm{U})=F_X^{-1}(l_X\cap \mathrm{U})\subset \mathbb{R}$, which are open intervals and make a finite or countable set $(\mathrm{I}_j)$ \ \{Ref. \cite{DieudonneTome1}, \S (3.19.6)\}. \hfill $\square$ \\ 

\paragraph{An argument of transversality.}\label{Conjecture} Assume that $v$ does not vanish and that each maximal integral curve is non-periodic and is closed in $\mathrm{V}$. Given an arbitrary point $X\in \mathrm{V}$, let $\chi: \mathrm{U}_0 \rightarrow \mathrm{I}\times \Omega_0$ be a straightening-out chart in the neighborhood of $X$, and let $\mathrm{E}$ be the closed countable subset of $\Omega_0 $, having only isolated points, such that $\chi (l_X\cap \mathrm{U}_0) = \mathrm{I}\times \mathrm{E}$. As \hyperref[Proposition 3]{Proposition 3} states: by restricting $\chi $ to an open subset $\mathrm{U}\subset \mathrm{U}_0$ for which $\chi (\mathrm{U})$ has the form (\ref{chi(U)}) with  $\Omega \subset \Omega_0 $ any open neighborhood of ${\bf x}\equiv P_S(\chi (X))$ such that $\mathrm{E}\cap \Omega =\{{\bf x}\}$, we ensure that $l_X\cap \mathrm{U}$ is connected. From Point (iii) of \hyperref[Proposition 2]{Proposition 2}, it follows that we will obtain a nice $v$-adapted chart by restricting $\chi $ to an open neighborhood $\mathrm{W}\subset\mathrm{U} $ of $X$ --- if it exists --- such that, for any $Y\in \mathrm{W}$, $l_Y\cap \mathrm{U}$ be connected. To this purpose, we observe that, when the boundary of  $\mathrm{U}$ is a regular hypersurface, then the same should be true for the boundary of the open set $\mathcal{D}_\mathrm{U}\equiv F^{-1}(\mathrm{U})\subset \mathbb{R}\times \mathrm{V}$: that boundary $\Sigma _\mathrm{U}\equiv \mathrm{Fr}(\mathcal{D}_\mathrm{U})$ should ``normally" be a hypersurface in $\mathbb{R}\times \mathrm{V}$. I.e., $\Sigma _\mathrm{U}$ should be a submanifold of codimension 1 of the $(N+2)-$dimensional manifold $\mathbb{R}\times \mathrm{V}$. Then, ``generically", a straight line $\mathbb{R}\times \{X\}\subset \mathbb{R}\times \mathrm{V}$ that intersects $\Sigma _\mathrm{U}$ is nowhere tangent to it, i.e. it is transverse to it at each intersection point. Thus, when $Y$ is sufficiently close to $X$, the intersection points should be slightly displaced but should remain in the same number, hence the structure of $\mathrm{I}'_{Y \mathrm{U}}\times \{Y\}=(\mathbb{R}\times \{Y\})\cap \mathcal{D}_\mathrm{U}$ should be the same as for $X$. In particular, if $\mathrm{I}'_{X \mathrm{U}} $ is an interval, i.e. (according to \hyperref[Proposition 4]{Proposition 4}) if $l_X\cap \mathrm{U}$ is connected, then also $\mathrm{I}'_{Y \mathrm{U}}$ should be an interval, i.e., also $l_Y\cap \mathrm{U}$ should be connected. \\

To make this line of reasoning precise, we first introduce, for a general hypersurface $\Sigma $ in $\mathbb{R}\times \mathrm{V}$, the set $\mathrm{S}_\Sigma$ of the points $X$ in V for which the straight line $\mathbb{R}\times \{X\}$ is tangent to $\Sigma $ at one point at least. (Our aim is to use this when $\Sigma $ is the boundary hypersurface $\Sigma _\mathrm{U}\equiv \mathrm{Fr}(\mathcal{D}_\mathrm{U})$ introduced above.) Note that the tangent space to $\mathbb{R}\times \mathrm{V}$ at some point $(s,X)$ is $\mathrm{T}_{(s,X)}(\mathbb{R}\times \mathrm{V})\simeq \mathrm{T}_s \mathbb{R} \times \mathrm{T}_X \mathrm{V}\simeq \mathbb{R}\times \mathrm{T}_X \mathrm{V}$. The tangent space to $\mathbb{R}\times \{X\}$ at $(s,X)$ is $\mathbb{R}\xi $, where $\xi \equiv (1,0_X)$, $0_X$ being the zero element of $\mathrm{T}_X \mathrm{V}$. Thus we define
\be
\mathrm{S}_\Sigma  \equiv \{X \in \mathrm{V};\,\exists s \in \mathbb{R}: (s,X)\in \Sigma \ \mathrm{and}\ (1,0_X) \in \mathrm{T}_{(s,X)}\Sigma  \}.
\ee

\vspace{4mm}
However, we note that $\mathcal{D}_\mathrm{U}\equiv F^{-1}(\mathrm{U})$ is in general not a bounded domain of $\mathbb{R}\times \mathrm{V}$, even if the open set $\mathrm{U}\subset \mathrm{V}$ is bounded. 
\footnote{\label{ConstantVectorField}\ 
For instance, consider a constant vector field $v$ on $\mathrm{V}\equiv \mathbb{R}^n$. Then we have simply $F(s,X)=X+sv$. Hence, for any $X_0\in \mathrm{V}$, its inverse image is the unbounded straight line $F^{-1}(X_0)=\{(s,X)\in \mathbb{R}\times \mathrm{V};\, s\in \mathbb{R}, X=X_0-sv$\}. From this, one may deduce more general examples by applying a diffeomorphism.
}
Therefore, even if $X\in \mathrm{V}$ is such that the straight line $\mathbb{R}\times \{X\}$ is not tangent to the boundary hypersurface $\Sigma_\mathrm{U}$, i.e. $X \not \in \mathrm{S}_{\Sigma_\mathrm{U}} $, it may happen that $\mathbb{R}\times \{X\}$ is ``tangent to $\Sigma_\mathrm{U} $ at infinity", in which case any line $\mathbb{R}\times \{Y\}$, however close $Y$ can be to $X$, may have ``new" intersections with the hypersurface $\Sigma_\mathrm{U} $. For a general hypersurface $\Sigma $ in $\mathbb{R}\times \mathrm{V}$, in fact for {\it any subset} $\Sigma $ of $\mathbb{R}\times \mathrm{V}$, we thus introduce:
\be\label{Def S_infini}
\mathrm{S}_{\Sigma\,\infty}  \equiv \{X \in \mathrm{V}; \lim_{r\rightarrow \infty } \mathrm{inf}_{\,\abs{s}\ge r}\, d(X,\Sigma _s) =0\},
\ee
where, for any subset $\mathcal{B}$ of $\mathbb{R}\times \mathrm{V}$, we define $\mathcal{B}_s$ to be its slice in V at $s\in \mathbb{R}$:
\be
\mathcal{B}_s\equiv \{X\in \mathrm{V};\,(s,X)\in \mathcal{B}\}\subset \mathrm{V},
\ee
and where, given any distance $d$ that generates the (metrizable) topology of V, one defines for any subset A of V and any point $X\in \mathrm{V}$: $d(X,\mathrm{A})\equiv \inf_{Z\in \mathrm{A}}\, d(X,Z)$. 

\subsection{Relevant theorems}\label{Transverse}

\paragraph{Theorem 1.}\label{Theorem 1} {\it Let $\Sigma $ be a hypersurface of $\mathbb{R}\times \mathrm{V}$ that is closed in $\mathbb{R}\times \mathrm{V}$. Let $\mathrm{K}=[\alpha,\beta ]$ be a compact interval of $\mathbb{R}$ ($\alpha < \beta $). Suppose that, for some $X_0 \in \mathrm{V}$, the intersection $(\mathrm{K}\times \{X_0\})\cap \Sigma $ be a singleton $(s_0,X_0)$ with $\alpha <s_0<\beta $ and that this intersection be transverse, i.e. $X_0 \not \in \mathrm{S}_\Sigma $. Then there is a neighborhood $\mathrm{W}$ of $X_0$ such that, for any $X\in \mathrm{W}$, $(\mathrm{K}\times \{X\})\cap \Sigma $ is a singleton $\Phi(X)$ and this intersection is transverse. Moreover, the function $\Phi :\mathrm{W}\rightarrow \Sigma $ is smooth.} \\

\noi {\it Proof.} Since $(s_0,X_0)\in \Sigma $, and since $\Sigma $ is a submanifold of dimension $N+1$ of the $(N+2)-$dimensional differentiable manifold $\mathbb{R}\times \mathrm{V}$, there is a chart $(\Psi ,\mathcal{U})$ on $\mathbb{R}\times \mathrm{V}$, with $(s_0,X_0)\in \mathcal{U}$, for which the function 
\be
g: \mathcal{U}\rightarrow \mathbb{R},\ (s,X)\mapsto g(s,X)\equiv \Psi ^{N+2}(s,X)
\ee
[the last component of $\Psi (s,X)$] is such that 
\be\label{Sigma:g=0}
(s,X)\in \mathcal{U}\cap \Sigma\ \Longleftrightarrow \ g(s,X)=0
\ee
\{e.g. \cite{DieudonneTome3}, \S (16.8.3)\}. Then, if $(s,X)\in \mathcal{U}\cap\Sigma $, a vector $\eta  =(a,u)\in \mathrm{T}_{(s,X)}(\mathbb{R}\times \mathrm{V})\simeq \mathbb{R}\times \mathrm{T}_X \mathrm{V}$ is in the tangent space $\mathrm{T}_{(s,X)}\Sigma $, iff $\dd g_{(s,X)}(\eta )= 0$. Since the intersection $(\mathrm{K}\times \{X_0\})\cap \Sigma =\{(s_0,X_0)\}$ is transverse, we have
\be\label{dg(xi)ne0}
\dd g_{(s_0,X_0)}(\xi _0)\ne 0, \quad \xi _0\equiv (1,0_{\mathrm{T}_{X_0}\mathrm{V}}). 
\ee
Since $\mathrm{K}$ is a neighborhood of $s_0$: replacing $\mathcal{U}$ with a smaller neighborhood of $(s_0,X_0)$ if necessary, we may assume that, if $(s,X)\in \mathcal{U}$, then $s\in \mathrm{K}$. Consider a chart $(\chi ,\mathrm{U})$ of V in a neighborhood U of $X_0$, thus $\chi (X)={\bf X}=(x^\mu )\in \mathbb{R}^{N+1}$ for $X\in \mathrm{U}$. The mapping $\Xi :(s,X)\mapsto (s,\chi (X))$ is a local chart of $\mathbb{R}\times \mathrm{V}$ defined in the neighborhood $\mathbb{R}\times \mathrm{U}$ of $(s_0,X_0)$, while $\Psi$ is also a local chart of $\mathbb{R}\times \mathrm{V}$, defined in the neighborhood $\mathcal{U}$ of $(s_0,X_0)$. Let 
\be\label{f local g in Ixchi}
f(s,{\bf X})\equiv g(\Xi ^{-1}(s,{\bf X}))=g(s,\chi ^{-1}({\bf X}))
\ee
be the local expression of $g$ in the chart $\Xi $. This function $f$ is defined and $\mathcal{C}^\infty $ on the open subset $\mathcal{O}\equiv \Xi ((\mathbb{R}\times \mathrm{U})\cap \mathcal{U})$ of $\mathbb{R}^{N+2}$. Note that $\mathcal{O}$ contains $(s_0,{\bf X}_0)$, where ${\bf X}_0\equiv \chi (X_0)$. With this local expression, we have for any vector $\eta  =(a,u)\in \mathrm{T}_{(s,X)}(\mathbb{R}\times \mathrm{V})$:
\be\label{local dg}
\dd g_{(s,X)}(\eta )=\frac{\partial f}{\partial s}(s,{\bf X})\,a+\frac{\partial f}{\partial x^\mu }(s,{\bf X})\,u ^\mu ,
\ee
where ${\bf X}\equiv \chi (X)$, and where $(a,(u^\mu ))\ (\mu =0,...,N)$ are the components of $\eta $ in the product chart $\Xi $. From (\ref{Sigma:g=0}), we have 
\be
f(s_0,{\bf X}_0)=0.
\ee
From (\ref{dg(xi)ne0}) and (\ref{local dg}) we have 
\be\label{transverse-inter-with-f}
\dd g_{(s_0,X_0)}(\xi _0)=\frac{\partial f}{\partial s}(s_0,{\bf X}_0) \ne 0.
\ee
We can thus apply the implicit function theorem: there is an open neighborhood $\mathrm{A}$ of ${\bf X}_0$ and a unique smooth function $\varphi : \mathrm{A}\rightarrow \mathbb{R}$, such that $\varphi ({\bf X}_0)=s_0$ and that, for any ${\bf X}\in \mathrm{A}$:
\be\label{implicit f}
(\varphi ({\bf X}),{\bf X})\in \mathcal{O},\ f(\varphi ({\bf X}),{\bf X})=0, \ \frac{\partial f}{\partial s}(\varphi ({\bf X}),{\bf X}) \ne 0.
\ee
We define a smooth mapping $\Phi $ by setting for any $X\in \mathrm{W}'\equiv \chi ^{-1}(\mathrm{A})\subset \mathrm{V}$: 
\be\label{Def Phi}
\Phi (X)\equiv (\varphi (\chi (X)),X).
\ee
Thus, for any $X\in \mathrm{W}'$, with ${\bf X}\equiv \chi (X)\in \mathrm{A}$, we have $\Phi (X)=(\varphi ({\bf X}),\chi ^{-1}({\bf X}))\\=\Xi ^{-1}(\varphi ({\bf X}),{\bf X})\in \Xi ^{-1}(\mathcal{O})\subset \mathcal{U}$. In particular, $\Phi (X_0)=(s_0,X_0)$. For any $X\in \mathrm{W}'$, we have from (\ref{f local g in Ixchi}), (\ref{implicit f}) and (\ref{Def Phi}): $g(\Phi (X))=0$, thus  $\Phi (X)\in \Sigma \cap \mathcal{U}$. As with (\ref{dg(xi)ne0}) and (\ref{transverse-inter-with-f}), (\ref{implicit f}) means that the intersection $\Phi (X)\in (\mathrm{K}\times \{X\})\cap \Sigma$ is transverse. The definition (\ref{Def Phi}) entails also that $\mathrm{Pr}_2 \circ \Phi = \mathrm{Id}_{\mathrm{W}'}$, where $\mathrm{Pr}_2 : \mathbb{R} \times \mathrm{V}\rightarrow \mathrm{V},\ (s,X)\mapsto X$. Thus, the rank of $\mathrm{Pr}_2 \circ \Phi$ is $\mathrm{dim}\,\mathrm{V}=N+1$, and since the rank of $\mathrm{Pr}_2 \circ \Phi $ is not larger than that of $\Phi $, this latter is also $N+1=\mathrm{dim}\,\Sigma $, i.e. $\Phi $ is a submersion. (All of this is true at any point $X\in \mathrm{W}'$.) It follows that $\mathrm{W}''\equiv \Phi(\mathrm{W}')$ is an open set in the manifold $\Sigma $ \{\cite{DieudonneTome3}, \S (16.7.5)\}, hence is a neighborhood of $(s_0,X_0)=\Phi (X_0)$ in $\Sigma $. We have
\be\label{(s,X)=Phi(X)}
\forall (s,X)\in \mathrm{W}'',\ (s,X)=\Phi (X).
\ee
We claim that there is some neighborhood $\mathrm{W}\subset \mathrm{W}'$ of $X_0$, such that
\be\label{KxX inter Sigma = Phi(X)}
\forall X\in \mathrm{W},\ [s\in \mathrm{K}\ \mathrm{and} \ (s,X)\in \Sigma ]\Longrightarrow [(s,X)=\Phi (X)]. 
\ee
Note that, if $\mathrm{W}\subset \mathrm{W}'$, the reverse implication is also true for any $X\in \mathrm{W}$, by the construction of the mapping $\Phi $. Thus, if (\ref{KxX inter Sigma = Phi(X)}) is true, $\mathrm{W}$ is as stated by \hyperref[Theorem 1]{Theorem 1}. We will reason {\it ab absurdo}. If it does not exist such a neighborhood $\mathrm{W}$, then, for any integer $n>0$, one can find $s_n\in \mathrm{K}$ and $X_n\in \mathrm{W}'\cap \mathrm{B}(X_0,1/n)$ such that $(s_n,X_n)\in \Sigma $ and $(s_n,X_n) \ne \Phi (X_n)$. Thus $X_n\rightarrow X_0$ and, by extraction in the compact set $\mathrm{K}$, we may assume that $s_n$ has a limit $s\in \mathrm{K}$, so that $(s_n,X_n)\rightarrow (s,X_0)$. If it happened that $s=s_0$, then, since $\mathrm{W}''$ is a neighborhood of $(s_0,X_0)$ and $(s_n,X_n)\rightarrow (s,X_0)$, we would have $(s_n,X_n)\in \mathrm{W}''$ for large enough $n$, hence $(s_n,X_n)=\Phi (X_n)$ from (\ref{(s,X)=Phi(X)}), which is a contradiction. Thus $s\ne s_0$. But since $\Sigma $ is closed, we have $(s,X_0)\in \Sigma $, which contradicts the assumption that $(\mathrm{K}\times \{X_0\})\cap \Sigma =\{(s_0,X_0)\}$. This completes the proof of \hyperref[Theorem 1]{Theorem 1}.  \hfill $\square$ 

\paragraph{Theorem 2.}\label{Theorem 2} {\it Let $\mathcal{U}$ be an open domain in $\mathbb{R}\times \mathrm{V}$, whose boundary $\Sigma $ be a  hypersurface of $\mathbb{R}\times \mathrm{V}$. Assume that, for some $X\in \mathrm{V}$, all intersections of the straight line $\mathcal{L}\equiv  \mathbb{R}\times \{X\}$ with $\Sigma $ are transverse. Then the boundary of $ \mathcal{L}\cap \mathcal{U}$ in $ \mathcal{L}$ is $\mathrm{Fr}_{(\mathcal{L})}(\mathcal{L} \cap \mathcal{U})=\mathcal{L} \cap \Sigma $.}\\

\noi {\it Proof.} The boundary of a subset is of course relative to which containing set is considered. Here, the boundary $\Sigma $ of the open set $\mathcal{U}\subset \mathbb{R}\times \mathrm{V}$ is relative to the whole manifold $\mathcal{V}\equiv \mathbb{R}\times \mathrm{V}$, i.e. $\Sigma \equiv  \overline{\mathcal{U}}\cap \overline{\complement\mathcal{U}}$, where the upper bar $\overline{.}$ denotes the adherence in $\mathcal{V}$ and $\complement$ the complementary set in $\mathcal{V}$. Thus 
\be\label{Sigma = D bar - D}
\Sigma =  \overline{\mathcal{U}}\setminus  \mathcal{U},
\ee 
since $\mathcal{U}$ is open. The boundary of some subset $ \mathcal{B}\subset  \mathcal{L}$ in $ \mathcal{L}$ (or relative to $ \mathcal{L}$) is $\mathrm{Fr}_{(\mathcal{L})}(\mathcal{B}) \equiv \overline{\mathcal{B}}^\mathcal{L}\cap \overline{\complement_\mathcal{L}\mathcal{B}}^\mathcal{L} $, where $\overline{\mathcal{B}}^\mathcal{L}$ denotes the adherence of $\mathcal{B}$ in $\mathcal{L}$, and where $\complement_\mathcal{L}\mathcal{B}\equiv \mathcal{L}\setminus \mathcal{B}$ is the complementary set of $\mathcal{B}$ in $\mathcal{L}$. However, here $\mathcal{L}\equiv  \mathbb{R}\times \{X\}$ is closed in $\mathcal{V}=\mathbb{R}\times \mathrm{V}$, hence we have $\overline{\mathcal{B}}^\mathcal{L}=\overline{\mathcal{B}}$, the adherence in the whole set $\mathcal{V}$. Thus
\be\label{Fr_A A inter D}
\mathrm{Fr}_{(\mathcal{L})}(\mathcal{L} \cap \mathcal{U})   \equiv  \overline{\mathcal{L} \cap \mathcal{U}\,}^\mathcal{L} \cap \overline{\complement_\mathcal{L}\left( \mathcal{L} \cap \mathcal{U}\right)\,}^\mathcal{L}=
\overline{\mathcal{L} \cap \mathcal{U}} \cap \overline{\complement_\mathcal{L}\mathcal{U}}.
\ee

\vspace{2mm}
\noi Again because $\mathcal{L}$ is closed in $\mathcal{V}$, we have  
$\overline{\mathcal{L} \cap \mathcal{U}} \subset \mathcal{L} \cap 
\overline{\mathcal{U}}$. We will show that we have exactly
\be\label{(A inter D) bar = A inter D bar}
\overline{\mathcal{L} \cap \mathcal{U}} = \mathcal{L} \cap 
\overline{\mathcal{U}}.
\ee
We shall in fact show that 
\be\label{A inter Sigma inclus dans (A inter D)bar}
\mathcal{L}\cap \Sigma  \subset \overline{\mathcal{L} \cap \mathcal{U}}.
\ee
Since (\ref{Sigma = D bar - D}) implies that $\left(\mathcal{L} \cap \overline{\mathcal{U}}\right)\setminus \overline{\mathcal{L} \cap \mathcal{U}} \subset \left(\mathcal{L} \cap \overline{\mathcal{U}}\right)\setminus \left(\mathcal{L} \cap \mathcal{U} \right)=\mathcal{L}\cap \Sigma $, this will prove that $\mathcal{L} \cap 
\overline{\mathcal{U}} \subset \overline{\mathcal{L} \cap \mathcal{U}} $, whence (\ref{(A inter D) bar = A inter D bar}).\\

To prove (\ref{A inter Sigma inclus dans (A inter D)bar}), consider an arbitrary point $p_0=(s_0,X)\in \mathcal{L}\cap \Sigma$. As in the proof of \hyperref[Theorem 1]{Theorem 1}, let $(\Psi ,\mathcal{W})$ be a chart on $\mathcal{V}$, with $p_0\in \mathcal{W}$, such that $p\in \mathcal{W}\cap \Sigma\Leftrightarrow g(p)=0,$ where $g= \Psi ^{n}$, with $n=\mathrm{dim}(\mathcal{V})=N+2$. Since we assume that the intersection $p_0\in \mathcal{L}\cap \Sigma $ is transverse, we have again (\ref{dg(xi)ne0}). Hence, there is an interval $\mathrm{J}=]s_0-r,s_0+r[$ in which $s=s_0$ is the only zero of the smooth function $\varphi (s)\equiv g(s,X)$. (That is, $p_0$ is the only intersection of $\mathrm{J}\times \{X\}$ with $\Sigma $.) Thus, we may assume that, say, $g(s,X)>0 \mathrm{\ for\ }s_0<s<s_0+r$, so that
\be\label{left in U_+}
g(s,X)<0 \mathrm{\ for\ }s_0-r<s<s_0.
\ee
Replacing $\Psi (p)$ by $\Psi (p)-\Psi (p_0)$, we may assume that $\Psi (p_0)={\bf 0}_ {\mathbb{R}^{n}}$. There is some open ball $\mathrm{W}=\mathrm{B}({\bf 0},r)\subset \Psi(\mathcal{W})$. Replacing $\mathcal{W}$ by $\Psi ^{-1}(\mathrm{W})$, we have that $\mathcal{W}_+\equiv \{ p\in \mathcal{W};x^n\equiv \Psi ^n(p)>0\}$ is just $\mathcal{W}_+=\Psi ^{-1}(\mathrm{W}_+)$, where $\mathrm{W}_+\equiv \{{\bf P}\in \mathrm{W};x^n>0\}$, hence $\mathcal{W}_+$ is non-empty and connected as is $\mathrm{W}_+$. The same is true for $\mathcal{W}_-\equiv \{ p\in \mathcal{W};\Psi^n(p)<0\}$. Since $p_0\in \Sigma \subset \overline{\mathcal{U}}$, and since $\mathcal{W}$ is a neighborhood of $p_0$, we have $\mathcal{U}\cap \mathcal{W} \ne \emptyset$, so let $p'_0\in \mathcal{U}\cap \mathcal{W}$. Because $\mathcal{W}$ is the disjointed union $\mathcal{W}=\mathcal{W}_+\cup \mathcal{W}_-\cup (\mathcal{W}\cap \Sigma) $, and because $p'_0\in \mathcal{U}$ cannot belong to $\Sigma =\overline{\mathcal{U}}\setminus  \mathcal{U}$, we have either $p'_0\in \mathcal{W}_+$ or $p'_0\in \mathcal{W}_-$. Let us assume that, for instance, $p'_0\in \mathcal{W}_-$, so that $\mathcal{U}\cap \mathcal{W}_- \ne \emptyset$. It follows that $\mathcal{W}_-\subset \mathcal{U}\cap \mathcal{W}$, for otherwise the connected set $\mathcal{W}_-$ would intersect both $\mathcal{U}$ and $\complement \mathcal{U}$, hence would intersect the boundary $\Sigma $ --- which is impossible, since we have $x^n=0$ in $\Sigma $, not $x^n<0$. Therefore, we get from (\ref{left in U_+}) that $]s_0-r,s_0[\times \{X\}\subset \mathcal{U}\cap \mathcal{W}$, so that $p_0\equiv (s_0,X)\in \overline{\mathcal{L} \cap \mathcal{U}}$. This proves (\ref{A inter Sigma inclus dans (A inter D)bar}).\\

Combining (\ref{Fr_A A inter D}) and (\ref{(A inter D) bar = A inter D bar}) gives us: 
\be\label{Fr_A A inter D-2}
\mathrm{Fr}_{(\mathcal{L})}(\mathcal{L} \cap \mathcal{U})   =\mathcal{L} \cap \overline{\mathcal{U}} \cap \overline{\complement_\mathcal{L}\mathcal{U}}.
\ee
But $\complement_\mathcal{L}\mathcal{U}=\mathcal{L}\cap \complement\mathcal{U}$ is closed in $\mathcal{V}$, for $\mathcal{L}$ is closed and $\mathcal{U}$ is open. Hence $\mathcal{L} \cap  \overline{\complement_\mathcal{L}\mathcal{U}} = \mathcal{L} \cap \complement_\mathcal{L}\mathcal{U}=\mathcal{L} \cap \complement\mathcal{U}=\mathcal{L} \cap \overline{\complement\mathcal{U}}$. Therefore, (\ref{Fr_A A inter D-2}) rewrites as
\be\label{Fr_A A inter D-3}
\mathrm{Fr}_{(\mathcal{L})}(\mathcal{L} \cap \mathcal{U})   =\mathcal{L} \cap \overline{\mathcal{U}} \cap \overline{\complement\mathcal{U}}\equiv \mathcal{L} \cap \Sigma ,
\ee
which proves \hyperref[Theorem 2]{Theorem 2}. \hfill $\square$ 

\paragraph{Remark 2.1.}\label{Remark 2.1} With small and straightforward modifications, the foregoing proof shows the result in the much more general case that $\mathbb{R}\times \mathrm{V}$ (with $\mathrm{V}$ a differentiable manifold) is replaced by a general differentiable manifold $\mathcal{V}$ and $\mathcal{L}$ is the range $\mathcal{L}=C(\mathrm{J})$, assumed closed, of a smooth curve $C: \mathrm{J} \rightarrow \mathcal{V}$ (with $\mathrm{J}$ an interval of $\mathbb{R}$), such that all intersections $p\in \mathcal{L}\cap \mathcal{V}$ are transverse. It is easy to see that the latter assumption is necessary.

\paragraph{Remark 2.2.}\label{Remark 2.2} In the course of the proof, the following intuitively obvious result was proved: Suppose that the line $\mathcal{L}=C(\mathrm{J})$ intersects transversely at point $p=C(s_0)$ the boundary $\Sigma $, assumed to be a regular hypersurface, of some open domain $\mathcal{U}\subset \mathcal{V}$. Then, among the two parts of $\mathcal{L}$: $s<s_0$ and $s>s_0$, at least one is such that, when $s$ is close enough to $s_0$, we have $C(s)\in \mathcal{U}$.

\paragraph{Theorem 3.}\label{Theorem 3} {\it Let $\mathrm{U}$ be an open subset of $\mathrm{V}$.} (i) {\it Assume that $\overline{F^{-1}(\mathrm{U})} \subset \mathcal{D}$, with $\mathcal{D}$ the domain of the flow $F$. (This is true, in particular, if the flow is complete, i.e. $\mathcal{D}=\mathcal{V}$.) Then we have:}
\be\label{Sigma = F-1(S)}
\mathrm{Fr}(F^{-1}(\mathrm{U}))=F^{-1}(\mathrm{Fr}(\mathrm{U})).
\ee
(ii) {\it If $\mathrm{Fr}(\mathrm{U})$ is a hypersurface of $\mathrm{V}$, then $F^{-1}(\mathrm{Fr}(\mathrm{U}))$ is a hypersurface of $\mathcal{V}\equiv \mathbb{R}\times \mathrm{V}$.}\\

\noi {\it Proof of Point} (i). First, consider the more general context that $\mathcal{V}$ and $\mathrm{V}$ are merely metric spaces and $F:\mathcal{D}\rightarrow \mathrm{V}$ is merely a continuous map, with $\mathcal{D}$ an open subset of $\mathcal{V}$. Then, if $\overline{F^{-1}(\mathrm{U})} \subset \mathcal{D}$, we have
\be\label{Sigma inclus dans F-1(S)}
 \mathrm{Fr}(F^{-1}(\mathrm{U}))\subset F^{-1}(\mathrm{Fr}(\mathrm{U})).
\ee
Indeed, since $\mathrm{U}$ is open in $\mathrm{V}$, we have $\mathrm{Fr}(\mathrm{U})=\overline{\mathrm{U}}^\mathrm{V}\setminus \mathrm{U}$, as with Eq. (\ref{Sigma = D bar - D}). Since $F$ is continuous on $\mathcal{D}$, $F^{-1}(\mathrm{U})$ is open in $\mathcal{D}$; hence, $\mathcal{D}$ being an open subset of $\mathcal{V}$, $F^{-1}(\mathrm{U})$ is open in $\mathcal{V}$. Therefore, we have similarly $\mathrm{Fr}(F^{-1}(\mathrm{U})) =\overline{F^{-1}(\mathrm{U})} \setminus F^{-1}(\mathrm{U})$. Again because $F$ is continuous on $\mathcal{D}$, we have $\overline{F^{-1}(\mathrm{U})}^\mathcal{D}\subset F^{-1}\left(\overline{\mathrm{U}}^\mathrm{V}\right)$. But $\overline{F^{-1}(\mathrm{U})}^\mathcal{D}=\overline{F^{-1}(\mathrm{U})}$ since  $\overline{F^{-1}(\mathrm{U})} \subset \mathcal{D} $. Thus $\mathrm{Fr}(F^{-1}(\mathrm{U})) \subset F^{-1}\left(\overline{\mathrm{U}}^\mathrm{V}\right)\setminus F^{-1}(\mathrm{U})$, whence (\ref{Sigma inclus dans F-1(S)}).\\

To prove the reverse inclusion, we consider any $p=(s,X)\in F^{-1}(\mathrm{Fr}(\mathrm{U}))$, thus $Y\equiv F(s,X)\in \mathrm{Fr}(\mathrm{U})$, and we will show that $p\in \mathrm{Fr}(F^{-1}(\mathrm{U}))$. There exists an open neighborhood $\mathcal{U}$ of $p$, having the form $\mathcal{U}=\mathrm{J}\times \mathrm{W}$, with $\mathrm{J}$ an open interval containing $s$ and $0$, and with $\mathrm{W}$ an open neighborhood of $X$ in $\mathrm{V}$, such that $\mathcal{U}\subset \mathcal{D}$ \{Ref. \cite{DieudonneTome4}, \S (18.2.5)\}. For all $t\in \mathrm{J}$, $F_t\equiv F(t,.)$ is a homeomorphism of $\mathrm{W}$ onto $\mathrm{W}_t\equiv F_t(\mathrm{W})$, moreover $F_{-t}$ is defined over $\mathrm{W}_t$ and is the inverse homeomorphism of $F_t$ [these two points result from (\ref{C_X vs C_Y})] \{Ref. \cite{DieudonneTome4}, \S (18.2.8)\}. (We note for the proof of Point (ii) that $F_t$, as well as $F_{-t}$, is of class $\mathcal{C}^q$ if the vector field $v$ is itself $\mathcal{C}^q$ \cite{DieudonneTome4}.) Now we consider any neighborhood $\mathcal{U}'$ of $p$ and we show that it intersects both $F^{-1}(\mathrm{U})$ and $\complement F^{-1}(\mathrm{U})$. We may assume that $\mathcal{U}'\subset \mathcal{U}$ and has the form $\mathcal{U}'=\mathrm{J}'\times \mathrm{W}'$, with $\mathrm{J}'\subset \mathrm{J}$ an open interval containing $s$, and with $\mathrm{W}'\subset \mathrm{W}$ an open neighborhood of $X$ in $\mathrm{V}$. Thus $F_s(\mathrm{W}')$ is an open neighborhood of $Y=F_s(X)$. Since $Y \in \mathrm{Fr}(\mathrm{U})$, both $F_s(\mathrm{W}')\cap \mathrm{U}$ and $F_s(\mathrm{W}')\cap \complement \mathrm{U}$ are non-empty, so let $Y'\in F_s(\mathrm{W}')\cap \mathrm{U}$ and $Y''\in F_s(\mathrm{W}')\cap \complement \mathrm{U}$. Therefore, $X'\equiv F_{-s}(Y')\in \mathrm{W}'$. Thus, $F(s,X')=Y'\in \mathrm{U}$, and also $(s,X')\in \mathrm{J}'\times \mathrm{W}'$, so that $(s,X')\in (\mathrm{J}'\times \mathrm{W}')\cap F^{-1}(\mathrm{U})$. In just the same way, with $X''\equiv F_{-s}(Y'')$, we get that $(s,X'')\in (\mathrm{J}'\times \mathrm{W}')\cap F^{-1}(\complement \mathrm{U})$. Since $F^{-1}(\complement \mathrm{U})=\complement_\mathcal{D} F^{-1}(\mathrm{U})\subset \complement F^{-1}(\mathrm{U})\equiv \mathcal{V} \setminus F^{-1}(\mathrm{U})$, we have shown that any neighborhood $\mathcal{U}'$ of $p$ intersects both $F^{-1}(\mathrm{U})$ and $\complement F^{-1}(\mathrm{U})$, thus $p\in \mathrm{Fr} (F^{-1}(\mathrm{U}))$.  Therefore, we have indeed $F^{-1}(\mathrm{Fr}(\mathrm{U}))\subset \mathrm{Fr} (F^{-1}(\mathrm{U}))$. Together with (\ref{Sigma inclus dans F-1(S)}), this proves (\ref{Sigma = F-1(S)}). \\ 

\noi {\it Proof of Point} (ii). Let $p=(s,X)\in F^{-1}(\mathrm{Fr}(\mathrm{U}))\subset \mathcal{V}$. Define  $\mathrm{J}$ (with $s\in \mathrm{J}$), $ \mathrm{W}\subset \mathrm{V}$ (with $X\in\mathrm{W}$), $\mathcal{U}=\mathrm{J}\times \mathrm{W}$, and $F_t\ (t\in \mathrm{J})$ just as at the beginning of the foregoing paragraph. Since $\mathrm{S}\equiv \mathrm{Fr}(\mathrm{U})$ is assumed to be a hypersurface of $\mathrm{V}$, and since by hypothesis $Y\equiv F(s,X)\in \mathrm{S}$, let $\chi: Y'\mapsto {\bf Y}\equiv (y^1,...,y^{n-1})$ be a chart of $\mathrm{V}$ (with $n-1=\mathrm{dim}(\mathcal{V})-1=\mathrm{dim}(\mathrm{V})$), defined in the neighborhood of $Y$, and such that 
\be\label{S: y^{n-1}=0}
Y'\in \mathrm{S}\cap \mathrm{Dom}(\chi)  \Leftrightarrow y^{n-1}\equiv \chi ^{n-1}(Y')=0. 
\ee
We may assume that $\mathrm{Dom}(\chi)$, the domain of $\chi $, is just $ \mathrm{W}_s\equiv F_s( \mathrm{W})$. Then the mapping 
\be
\Xi :\mathcal{U}=\mathrm{J}\times \mathrm{W}\rightarrow \mathbb{R}^n,\quad (t,X')\mapsto (t,\chi (F_s(X')))=(t,{\bf Y})
\ee
is a chart of $\mathcal{V}$ in the neighborhood of $p$. Consider the mapping 
\be\label{Psi}
\Psi :\mathcal{U}'\equiv F^{-1}(\mathrm{W}_s)\rightarrow \mathbb{R}^n,\quad (t,X')\mapsto (t,\chi (F(t,X'))).
\ee
(Note that $\mathcal{U}'$ is an open neighborhood of $p$.) We have (cf. Eq. (\ref{C_X vs C_Y})):
\be
F(t,X')=F(s+u,X')=F(u,F(s,X'))=F(u,F_s(X')),\quad u\equiv t-s.
\ee
Hence, the local expression of $\Psi $ in the chart $\Xi $ is:
\be\label{Def G}
 {\bf G}(t,{\bf Y})\equiv \Psi (\Xi ^{-1}(t,{\bf Y}))= (t,{\bf Z}(t,{\bf Y}))
\ee
with
\be
{\bf Z}(t,{\bf Y})=(z^1,...,z^{n-1})\equiv \chi (F(u,\chi ^{-1}({\bf Y}))).
\ee
Let ${\bf v}=(v^1,...,v^{n-1})={\bf v}({\bf Y})$ be the local expression of $v$ in the chart $\chi $. Thus ${\bf Z}$ is the value at $u\equiv t-s$ of the solution of $\frac{\dd {\bf Z}'}{\dd w}={\bf v}({\bf Z}'(w)),\ {\bf Z}'(w=0)= {\bf Y}$. 
Hence we have, uniformly w.r.t.  $(t,{\bf Y})$ in some neighborhood of $(s,{\bf Y}_0)\equiv \Xi (p)$:
\be
{\bf Z}(t,{\bf Y})={\bf Y}+(t-s){\bf v}({\bf Y})+O(u^2).
\ee
The Jacobian matrix of ${\bf G}(t,{\bf Y})=(t,{\bf Z}(t,{\bf Y}))$ at point $(s,{\bf Y}_0)$ is therefore:
\be\label{J}
J = \begin{pmatrix} 
1                  & 0       & ...    &  ...   & 0\\
v^1({\bf Y}_0)     & 1       & 0      & ...    & 0\\
...                & ...     & ...    & ...    & ...\\
v^{n-1}({\bf Y}_0) & 0       & 0      & ...    & 1\\
\end{pmatrix},
\ee
a triangular matrix with 1 on the diagonal, so $\det J=1$. It follows that $\Psi $ also is a chart of $\mathcal{V}$ in some neighborhood $\mathcal{U}''\subset \mathcal{U}\cap \mathcal{U}'$ of $p$. When $(t,X')\in \mathcal{U}''$, we have from (\ref{S: y^{n-1}=0}) and (\ref{Psi}):
\be\label{F^-1(S): Psi^{n-1}=0}
(t,X')\in F^{-1}(\mathrm{S})  \Leftrightarrow z^{n-1}\equiv \Psi ^{n}(t,X')\equiv \chi ^{n-1}(F(t,X'))=0. 
\ee
Hence, $F^{-1}(\mathrm{S})$ is a hypersurface of $\mathcal{V}$. The proof of \hyperref[Theorem 3]{Theorem 3} is complete.  \hfill $\square$ 

\paragraph{Remark 3.1.}\label{Remark 3.1} From (\ref{Psi}) and (\ref{Def G}), we have also ${\bf Z}(t,{\bf Y})=\chi (F(\Xi ^{-1}(t,{\bf Y})))$, thus ${\bf Z}(t,{\bf Y})$ is the local expression of $F$ in the charts $\Xi $ on $\mathcal{V}$ and $\chi $ on $\mathrm{V}$. Equation (\ref{J}) shows also that the Jacobian matrix of ${\bf Z}$ at $(s,{\bf Y}_0)\equiv \Xi (p)$ has rank $n-1=\mathrm{dim\,V}$. [Relation (\ref{S: y^{n-1}=0}), hence the fact that $p \in F^{-1}(\mathrm{S})$, are not used to get this; hence this is true at any point $p \in \mathcal{D}=\mathrm{Dom}(F)$.] Thus, $F$ is a submersion. Hence, it is transversal to any submanifold of $\mathrm{V}$. Point (ii) of \hyperref[Theorem 3]{Theorem 3} follows also from this \cite{GuilleminPollack-p28}.

\paragraph{Proposition 5.}\label{Proposition 5} {\it Let $\Sigma $ be a subset of $\mathbb{R}\times \mathrm{V}$. If $X\in \mathrm{V}\setminus \mathrm{S}_{\Sigma\,\infty}$ [cf. (\ref{Def S_infini})], there is a neighborhood $\mathrm{W}$ of $X$ and a real $R>0$ such that}
\be
\left(Y\in\mathrm{W} \mathrm{\ and\ }\abs{s}\ge R \right) \Rightarrow (s,Y) \not \in \Sigma.
\ee

\vspace{3mm}
\noi {\it Proof.} Clearly, $f(X)\equiv \lim_{r\rightarrow \infty } \mathrm{inf}_{\,\abs{s}\ge r}\, d(X,(\Sigma )_s)$ is well defined for any $X \in \mathrm{V}$ and verifies $0 \le f(X) \le +\infty $. Hence, if $X \not \in \mathrm{S}_{\Sigma\,\infty}$ it means that $f(X)>0$, or equivalently that there exists $R >0$ and $\delta >0$ such that
\be
\abs{s}\ge R \Rightarrow d(X,\Sigma _s) \ge \delta ,
\ee
hence $d(X,Z)\ge \delta$ is true for any $Z\in \Sigma _s$ if $\abs{s}\ge R$. We deduce from this that, if $Y\in \mathrm{V}$ verifies $d(X,Y) < \delta /2$, and if $\abs{s}\ge R$, we have $d(Y,Z)>\delta /2$ for any $Z\in \Sigma _s$, hence $d(Y,\Sigma _s)\ge \delta /2$. (Thus, $\mathrm{V}\setminus \mathrm{S}_{\Sigma\,\infty}$ is open in $\mathrm{V}$, in other words $\mathrm{S}_{\Sigma\,\infty}$ is closed.) In particular, if $d(X,Y) < \delta /2$ and $\abs{s}\ge R$, then $Y\not \in \Sigma _s$, i.e. $(s,Y) \not \in \Sigma$.  \hfill $\square$

\paragraph{Theorem 4.}\label{Theorem 4} {\it Let $\mathrm{U}$ be an open subset of $\mathrm{V}$ and set $\mathcal{D}_\mathrm{U}\equiv F^{-1}(\mathrm{U})$ and $\Sigma _\mathrm{U} \equiv \mathrm{Fr}(\mathcal{D}_\mathrm{U})$. Assume that $\Sigma _\mathrm{U} $ is a hypersurface of $\mathcal{V}\equiv \mathbb{R}\times \mathrm{V}$ which is closed in $\mathcal{V}$, that for some $X\in \mathrm{U}$, the open set $\mathrm{I}'_{X \mathrm{U}}$ is a bounded interval, and that $X \not \in (\mathrm{S}_{\Sigma _\mathrm{U} }\cup \mathrm{S}_{\Sigma _\mathrm{U} \,\infty }) $.}\\ (i) {\it There is a neighborhood $\mathrm{W}$ of $X$, $\mathrm{W}\subset \mathrm{U}$, such that for $Y\in \mathrm{W}$, also $\mathrm{I}'_{Y \mathrm{U}}$ is a bounded interval.} (ii) {\it If  $v$ does not vanish and all maximal integral curves are closed in $\mathrm{V}$ and non-periodic, then $l_Y\cap \mathrm{U}$ is connected for any $Y\in \mathrm{W}$.} (iii) {\it If moreover $\chi (\mathrm{U})$ has the form (\ref{chi(U)}) and there is a straightening-out chart $(\chi,\mathrm{U}_0) $ with $\mathrm{U}_0 \supset \mathrm{U}$, then the restriction of $\chi $ to $\mathrm{W}$ is a nice $v$-adapted chart.}\\

\noi {\it Proof.} (i) Since $X \not \in \mathrm{S}_{\Sigma _\mathrm{U} \,\infty } $, by \hyperref[Proposition 5]{Proposition 5} there is a neighborhood $\mathrm{W}_0$ of $X$ and a real $R>0$ such that
\be\label{W_0}
\left(Y\in\mathrm{W}_0 \mathrm{\ and\ }\abs{s}\ge R \right) \Rightarrow (s,Y) \not \in  \Sigma_\mathrm{U}.
\ee
Since $\mathrm{I}'_{X \mathrm{U}}$ is assumed to be a bounded interval (and $\mathrm{I}'_{X \mathrm{U}}\ne \emptyset$ since $0\in \mathrm{I}'_{X \mathrm{U}}$), let $\mathrm{I}'_{X \mathrm{U}}=]a,b[$, with $a,b\in \mathbb{R}, \,a<b$. In (\ref{W_0}) we may assume that $-R<a<b<R$. By \hyperref[Proposition 4]{Proposition 4}, we have 
\be
\mathrm{I}'_{X \mathrm{U}}\times\{X\}=]a,b[\times\{X\}=(\mathbb{R}\times \{X\})\cap \mathcal{D}_\mathrm{U}. 
\ee
Since by assumption $X \not \in \mathrm{S}_{\Sigma _\mathrm{U} } $, we get by \hyperref[Theorem 2]{Theorem 2} (setting $\mathcal{L}\equiv  \mathbb{R}\times \{X\}$):
\be
\{(a,X),(b,X)\}=\mathrm{Fr}_{(\mathcal{L})}\left(\mathcal{L}\cap \mathcal{D}_\mathrm{U}\right)=\mathcal{L}\cap\Sigma _\mathrm{U}.
\ee
Thus, considering the compact intervals $\mathrm{K}_1\equiv [-R,\frac{a+b}{2}]$ and $\mathrm{K}_2\equiv [\frac{a+b}{2},R]$, we have
\be\label{K_1 & K_2}
\left(\mathrm{K}_1\times \{X\}\right )\cap \Sigma _\mathrm{U}=\{(a,X)\},\quad \left(\mathrm{K}_2\times \{X\}\right )\cap \Sigma _\mathrm{U}=\{(b,X)\}.
\ee
Therefore, by \hyperref[Theorem 1]{Theorem 1}, there are neighborhoods $\mathrm{W_1}$ and $\mathrm{W_2}$ of $X$, and smooth functions $\Phi _j:\mathrm{W}_j\rightarrow \Sigma _\mathrm{U} \ (j=1,2)$, such that: 
\be\label{W_1 & W_2}
Y\in \mathrm{W}_j \Rightarrow (\mathrm{K}_j\times \{Y\})\cap \Sigma_\mathrm{U} =\{\Phi _j(Y)\}\quad (j=1,2),
\ee 
with transverse intersection. We have $\Phi _j(Y)=(\varphi _j(Y),Y)$, with $\varphi _j:\mathrm{W}_j\rightarrow \mathbb{R}$ a smooth function. From (\ref{K_1 & K_2}) and (\ref{W_1 & W_2}), we get: $\varphi _1(X)=a$ and $\varphi _2(X)=b$. Thus $\varphi _1(X)=a \ne b=\varphi _2(X)$; hence, by considering a small enough neighborhood $\mathrm{W}$ of $X$, with $\mathrm{W}\subset \mathrm{W}_0\cap \mathrm{W}_1\cap \mathrm{W}_2$, we get $\varphi _1(Y)\ne \varphi _2(Y)$ for $Y\in \mathrm{W}$. With (\ref{W_0}) and (\ref{W_1 & W_2}), this implies that 
\be\label{A_Y inter Sigma}
Y\in\mathrm{W} \Rightarrow (\mathbb{R}\times \{Y\})\cap \Sigma_\mathrm{U} = \{\Phi _1(Y),\Phi _2(Y)\}.
\ee 
Another application of \hyperref[Theorem 2]{Theorem 2} proves that 
\be
\mathrm{Fr}_{(\mathbb{R}\times \{Y\})}\left((\mathbb{R}\times \{Y\})\cap \mathcal{D}_\mathrm{U} \right)= (\mathbb{R}\times \{Y\})\cap \Sigma_\mathrm{U}.
\ee
Together with (\ref{A_Y inter Sigma}), this implies that $\mathrm{I}'_{Y \mathrm{U}}$, the open subset of $\mathbb{R}$ such that $\mathrm{I}'_{Y \mathrm{U}}\times\{Y\}=(\mathbb{R}\times \{Y\})\cap \mathcal{D}_\mathrm{U}$, has boundary $\{\varphi _1(Y),\varphi _2(Y)\}$. Therefore, $\mathrm{I}'_{Y \mathrm{U}}=]\varphi _1(Y),\varphi _2(Y)[$. \\
\noi (ii) This follows from Point (iii) in \hyperref[Proposition 4]{Proposition 4}.\\
\noi (iii) This follows from Point (iii) in \hyperref[Proposition 2]{Proposition 2}. \hfill $\square$ \\ 

\subsection{Adapted charts and ``normal" vector fields}

With \hyperref[Theorem 4]{Theorem 4}, we formalized our \hyperref[Conjecture]{transversality argument} in Subsect. \ref{Straight inter} to investigate the problem of the existence, in the neighborhood of any point $X\in \mathrm{V}$, of a nice $v$-adapted chart. Assuming that $v$ does not vanish and that all maximal integral curves are closed in $\mathrm{V}$ and non-periodic, let us check if \hyperref[Theorem 4]{Theorem 4} applies. Due to \hyperref[Proposition 4]{Proposition 4}, $\mathrm{I}'_{X \mathrm{U}}=F_X^{-1}(l_X\cap \mathrm{U})$ is an interval iff $l_X\cap \mathrm{U}$ is connected, and $\mathrm{I}'_{X \mathrm{U}}$ is bounded if $\mathrm{U}$ is relatively compact. As shown by \hyperref[Proposition 3]{Proposition 3}, the assumption that $l_X\cap \mathrm{U}$ is connected may be fulfilled by starting from a straightening-out chart $\chi: \mathrm{U}_0 \rightarrow \mathrm{I}\times \Omega_0$ in the neighborhood of the arbitrary point $X\in \mathrm{V}$ and by restricting $\chi $ to an open subset $\mathrm{U}\subset \mathrm{U}_0$ such that $\chi (\mathrm{U})$ has the form (\ref{chi(U)}) with  $\Omega \subset \Omega_0 $ a small enough open neighborhood of ${\bf x}\equiv P_S(\chi (X))$.\\

As shown by \hyperref[Theorem 3]{Theorem 3}, the assumption that $\Sigma _\mathrm{U} $ is a hypersurface of $\mathcal{V}$ that is closed in $\mathcal{V}$ is fulfilled, in particular, if the boundary of the open set $\mathrm{U}\subset \mathrm{V}$ is itself a hypersurface of $\mathrm{V}$ that is closed in $\mathrm{V}$, and if moreover $\overline{F^{-1}(\mathrm{U})} \subset \mathcal{D}$. The latter inclusion is true, in particular, if $\mathcal{D}=\mathcal{V}$, i.e., if the vector field $v$ is complete (in other words, if every maximal integral curve of $v$ is defined over the whole real line). Actually, this does not restrict in any way the set of the maximum integral curves, $\mathrm{N}_v\equiv \{l_X;\ X\in \mathrm{V}\}$ (the ``congruence of world lines", in the physical context with $N=3$). Indeed, there always exists a smooth function $\lambda : \mathrm{V}\rightarrow \mathbb{R}_+$, such that the vector field $\lambda v$ is complete, moreover the mappings $C_X$ corresponding to the maximal integral curves of $\lambda v$ are mere reparameterizations of those of $v$, so that the curves $l_X$ themselves are unchanged \cite{CalcutGompf2013}. Thus, the assumption that $\Sigma _\mathrm{U} $ is a hypersurface of $\mathcal{V}$ that is closed in $\mathcal{V}$ is not a restrictive one.\\

The assumption ``$X \not \in (\mathrm{S}_{\Sigma _\mathrm{U} }\cup 
\mathrm{S}_{\Sigma _\mathrm{U} \,\infty }) $" means that the straight line $\mathbb{R} \times \{X\}$ is not tangent to the hypersurface $\Sigma _\mathrm{U} $, and is not ``tangent to it at infinity". For a given hypersurface $\Sigma _\mathrm{U} $, the points thus excluded form a kind of apparent contour (of that hypersurface $\Sigma _\mathrm{U} $), having ``normally" measure zero in $\mathrm{V}$, hence this is true for a ``generic" point $X$. However, here the hypersurface $\Sigma _\mathrm{U}= F^{-1}(\mathrm{Fr(U)}) $ of $\mathbb{R}\times \mathrm{V}$ depends on the selected neighborhood $\mathrm{U}$ of the given point $X\in \mathrm{V}$. There is much freedom in the choice of this neighborhood, since it is merely required to have a regular boundary and have the form  (\ref{chi(U)}) with  $\Omega $ a small enough open neighborhood of ${\bf x}\equiv P_S(\chi (X))$. If it turns out that $X \in (\mathrm{S}_{\Sigma _\mathrm{U} }\cup \mathrm{S}_{\Sigma _\mathrm{U} \,\infty }) $ for some $\mathrm{U}$ satisfying these requirements, then a slightly deformed neighborhood $\mathrm{U}'$ does also satisfy them, but the boundary $\Sigma _{\mathrm{U}'}$ is also deformed w.r.t. $\Sigma _{\mathrm{U}}$. Hence, it seems plausible that, due to this freedom, there always exists $\mathrm{U}$ satisfying these requirements and such that $X \not \in (\mathrm{S}_{\Sigma _\mathrm{U} }\cup \mathrm{S}_{\Sigma _\mathrm{U} \,\infty }) $ --- unless $v$ has some ``pathology" that we were not able to describe in a more explicit way.\\

Thus, for any point $X\in \mathrm{V}$, the assumptions of \hyperref[Theorem 4]{Theorem 4} should be fulfilled in a suitable neighborhood $\mathrm{U}$ of $X$ if the vector field $v$ does not vanish, has all maximal integral curves closed in $\mathrm{V}$ and non-periodic, and does not suffer from the ``pathology" alluded to. Therefore, we set the following definition, the word ``normal" being justified by the foregoing discussion.

\paragraph{Definition 2.}\label{Normal} {\it A non-vanishing $\mathcal{C}^\infty $ vector field $\,v\,$ is called {\bf ``normal"} iff all maximal integral curves are closed in $\mathrm{V}$ and, moreover, any point $X \in \mathrm{V}$ has nested open neighborhoods $\mathrm{W}\subset \mathrm{U}\subset \mathrm{U}_0 $ such that:} (i) {\it There is a straightening-out chart $(\chi,\mathrm{U}_0) $.} (ii) {\it $\chi (\mathrm{U})$ has the form (\ref{chi(U)}).} (iii) {\it For any maximal integral curve $l$ of $v$ intersecting $\mathrm{W}$, the line $l\cap \mathrm{U}$ is connected.} \\

\noi The following result shows that this concept is relevant. Note that we do not need to assume that the maximal integral curves are non-periodic.

\paragraph{Theorem 5.}\label{Theorem 5} {\it Let $v$ be a non-vanishing global vector field, such that all maximal integral curves are closed in $\mathrm{V}$.} (i) {\it In order that, for any point $X\in \mathrm{V}$, there exist a nice $v$-adapted chart whose domain be an open neighborhood of $X$, it is necessary and sufficient that $v$ be normal.} (ii) {\it Also, in order that $v$ be normal, it is necessary and sufficient that any point $X \in \mathrm{V}$ have an open neighborhood $ \mathrm{W} $ on which there is a straightening-out chart $(\chi,\mathrm{W}) $, and such that, for any maximal integral curve $l$ of $v$, the line $l\cap \mathrm{W}$ is connected.}\\

\noi {\it Proof.} (i) The sufficiency is an immediate consequence of Point (iii) in \hyperref[Proposition 2]{Proposition 2}. Conversely, if for any $X\in \mathrm{V}$ there is a nice $v$-adapted chart $(\chi_1,\mathrm{U}_1) $ with $X \in \mathrm{U}_1$, then by Point (ii) of \hyperref[Proposition 0]{Proposition 0} we get a chart $(\chi,\mathrm{W}) $, with an open set $\mathrm{W}\subset \mathrm{U}_1$ and $X\in \mathrm{W}$, which (a) differs from $\chi_1$ merely by the time coordinate $y'^0$, and (b) is a straightening-out chart. From (a), $(\chi,\mathrm{W}) $ is also a $v$-adapted chart. Therefore, by Point (ii) of \hyperref[Proposition 2]{Proposition 2}: for any $Y\in \mathrm{W}$, the intersection $l_Y\cap \mathrm{W}$ is a connected set. This implies that for any maximal integral curve $l$ of $v$, the line $l\cap \mathrm{W}$ is a connected (possibly empty) set. By this and (b) above, and since $X\in \mathrm{V}$ is arbitrary, the vector field $v$ is normal. (This is the case $\mathrm{W}= \mathrm{U}= \mathrm{U}_0 $ in the definition.)\\

\noi (ii) For the reason just invoked, the condition is sufficient in order that $v$ be normal. Conversely, if $v$ is normal, consider any $X\in \mathrm{V}$. We know that there is a nice $v$-adapted chart $(\chi_1,\mathrm{U}_1) $ with $X \in \mathrm{U}_1$, and the proof of the necessity at Point (i) shows that from it we deduce a a straightening-out chart $(\chi,\mathrm{W}) $, with $X\in \mathrm{W}$, and such that for any maximal integral curve $l$ of $v$, the line $l\cap \mathrm{W}$ is a connected set.  \hfill $\square$\\

\paragraph{Examples.}\label{Examples normal} (i) Take $\mathrm{V}=\mathbb{R}^{N+1}$ and consider any constant vector field $v({\bf X})={\bf v}=\mathrm{Constant}\ne {\bf 0}$. The maximal integral curve at ${\bf X}\in \mathrm{V}$ is $l_{\bf X}=\{{\bf Y}={\bf X}+t{\bf v}; t\in \mathbb{R}\}$. To define a straightening-out chart explicitly, take ${\bf u}_1,...{\bf u}_N$ such that the vectors ${\bf u}_0\equiv {\bf v}, {\bf u}_1,...{\bf u}_N$ form a basis of $\mathbb{R}^{N+1}$, and define an invertible linear transformation $L$ of $\mathbb{R}^{N+1}$ by $L(x^\mu {\bf u}_\mu )=x^\mu {\bf e}_\mu$, where $({\bf e}_\mu )\ (\mu =0,...,N)$ is the canonical basis of $\mathbb{R}^{N+1}$. Now consider any open set of the form $\mathrm{U} =L ^{-1}(\mathrm{I}\times \Omega)$ with $\mathrm{I}=]-a,+a[$ and $\Omega $ an open subset of $\mathbb{R}^N$: we have explicitly $\mathrm{U} =\{s{\bf v}+x^j {\bf u}_j;s\in \mathrm{I},\ {\bf x}\equiv (x^j)\in \Omega \}$. The restriction $\chi $ of $L$ to $\mathrm{U}$ defines a straightening-out chart, because $L({\bf v})={\bf e}_0$ means that $v=\partial _0$ in that chart. Moreover, given ${\bf X}=s{\bf v}+x^j {\bf u}_j \in \mathrm{U}$ (thus $s\in \mathrm{I},\ {\bf x}\equiv (x^j) \in \Omega$), a point ${\bf Y}={\bf X}+t{\bf v}$ of $l_{\bf X}$ belongs to $\mathrm{U}$, iff $s+t \in \mathrm{I}$, so $l_{\bf X}\cap \mathrm{U}$ is connected.  Hence, a constant vector field is normal. And indeed, we have for  ${\bf Y}\in l_{\bf X}\cap \mathrm{U}$: $\chi ({\bf Y})=L({\bf X}+t{\bf v})=\chi ({\bf X})+t{\bf e}_0$, hence $\chi $ is $v$-adapted. Moreover, let $l\in \mathrm{D}_\mathrm{U}$, where $\mathrm{D}_\mathrm{U}$ is defined in Eq. (\ref{Def D_U}). Thus $l=l_{\bf X}$, where ${\bf X}=s{\bf v}+x^j {\bf u}_j \in \mathrm{U}$, i.e. $s\in \mathrm{I},\ {\bf x}\equiv (x^j)\in \Omega$. Then we have from the definition (\ref{Def chi bar}): $\bar{\chi }(l)={\bf x}=L({\bf x}_0)$ with ${\bf x}_0\equiv x^j {\bf u}_j$. Hence, the mapping $\bar{\chi }: \mathrm{D}_\mathrm{U}\rightarrow \mathbb{R}^N,\ l=l_{\bf X}=l_{s{\bf v}+{\bf x}_0}\mapsto {\bf x}=L({\bf x}_0)$ is injective, i.e., the $v$-adapted chart $\chi $ is nice.\\

(ii) If $\phi : \mathrm{V}\rightarrow \mathrm{V}'$ is a diffeomorphism, set $v'\equiv \phi ^*v:\ v'(X')\equiv D\phi_{\phi ^{-1}(X')} (v(\phi ^{-1}(X'))),\ X'\in \mathrm{V}'$. The maximal integral curves of $v'$, as well as the associate flow $F'$, are just the images of their counterparts for $v$: $F'(s,X')\equiv C'_{X'}(s)=\phi (C_{\phi ^{-1}(X')}(s))=\phi (F(s,\phi ^{-1}(X')))$. Moreover, if $(\chi,\mathrm{U}_0) $ is a straightening-out chart in the neighborhood of $X\in \mathrm{V}$, then so is $(\chi \circ \phi ^{-1},\phi (\mathrm{U}_0)) $ in the neighborhood of $X'\equiv \phi (X)\in \mathrm{V}'$. Therefore, if $v$ is normal, so is $v'$. \\

(iii) If we have a normal vector field $v$ on some differentiable manifold $\mathrm{V}$, having non-periodic orbits, and if $\mathrm{U}$ is any open subset of $\mathrm{V}$, let us show that the restriction $v'\equiv v_{\mid \mathrm{U}}$ is a normal vector field on the differentiable manifold $\mathrm{U}$. Given any $X\in \mathrm{U}$, it is easy to check from the definitions that the maximal open interval $\mathrm{I}'_X$ defining the orbit (maximal integral curve) $l'_X$ of $v'$ at $X$ is the connected component of $0$ in $\mathrm{I}'_{X \mathrm{U}}$ , the latter being defined in Eq. (\ref{Def I'_X}). (This result is true for any vector field.) It follows by Point (iii) in \hyperref[Proposition 4]{Proposition 4} that $l'_X$ is the connected component of $X$ in $l_X\cap \mathrm{U}$, where $l_X$ is the orbit of $v$ in $\mathrm{V}$ at $X$. Therefore, since the orbits of $v$ are closed subsets of $\mathrm{V}$, the orbits of $v'$ are closed subsets of $\mathrm{U}$. If $X \in \mathrm{U}$, there exists by Point (ii) of \hyperref[Theorem 5]{Theorem 5} a straightening-out chart $(\chi ,\mathrm{W})$ of $(\mathrm{V},v)$, with $X\in \mathrm{W}$ and $\chi (\mathrm{W})=\mathrm{I}\times \Omega $, such that for any $Y\in \mathrm{W}$, $l_Y\cap \mathrm{W}$ is connected. Let $\mathrm{B}=]x^0-r,x^0+r[\times ...\times ]x^N-r,x^N+r[$ be a ball centered at $\chi (X)=(x^\mu)$ and such that $\mathrm{B}\subset \chi (\mathrm{U})$, and set $\mathrm{W}'\equiv \chi ^{-1}(\mathrm{B})$. Then the restriction $\chi '\equiv \chi _{\mid \mathrm{W}'}$ is a straightening-out chart for $v'$ (up to a shift in $y^0$). For $Y\in \mathrm{W}'$, set ${\bf y}\equiv P_S(\chi (Y))$. Since $l_Y\cap \mathrm{W}$ is connected, we have $l_Y\cap \mathrm{W}=\chi ^{-1}(\mathrm{I}\times \{{\bf y}\})$ by Point (i) of \hyperref[Proposition 2]{Proposition 2}. Hence $l_Y\cap \mathrm{W}' =\chi ^{-1}(]x^0-r,x^0+r[\times \{{\bf y}\})$, thus a connected set $\subset l_Y\cap \mathrm{U}$. Since $l'_Y$ is the connected component of $Y$ in $l_Y\cap \mathrm{U}$, we have therefore $l'_Y\cap \mathrm{W}'=l_Y\cap \mathrm{W}'=\chi ^{-1}(]x^0-r,x^0+r[\times \{{\bf y}\})$. The conclusion follows by Point (ii) of \hyperref[Theorem 5]{Theorem 5}. \hfill $\square$\\

(iv) By combining the three former examples, we get that, if a manifold $\mathrm{V}$ is diffeomorphic to an {\it open subset} $\Gamma $ of $\mathbb{R}^{N+1}$ and $\phi :\Gamma \rightarrow \mathrm{V}$ is any diffeomorphism, then for any constant vector field ${\bf v}\ne {\bf 0}$ on $\Gamma $, its pushforward vector field by $\phi $, $v=\phi ^* {\bf v}$, is a normal vector field on $\mathrm{V}$. In the application to physics (for which $N=3$, as far as we know), this describes already a wide variety of spacetimes and vector fields, together with the associated reference fluids. Those are deformable in a very general way with respect to each other, by changing $\phi $, i.e. by transforming the integral curves by any diffeomorphism of $\mathrm{V}$. Of course, we expect that much more general normal vector fields do exist, due to the discussion at the beginning of this subsection. 

\section{The set of orbits of $v$ as a differentiable manifold}\label{N_v manifold}

The set $\mathrm{N}_v$ of the maximal integral curves of $v$ has been defined in Eq. (\ref{Def N_v}). In this section, we will show that, when $v$ is a normal vector field on the differentiable manifold $\mathrm{V}$, the set $\mathrm{N}_v$ can be endowed with a canonical structure of differentiable manifold.

\paragraph{Proposition 6.}\label{Proposition 6} {\it Let $v$ be a normal vector field on $\mathrm{V}$. Define the set $\mathcal{F}_v$ made of all nice $\,v$--adapted charts on $\mathrm{V}$. For any chart $\chi \in \mathcal{F}_v$, with domain $\mathrm{U}\subset \mathrm{V}$, let $\mathrm{D}_\mathrm{U}$ be defined by (\ref{Def D_U}), and, for any subset $\mathcal{O} \subset \mathrm{N}_v$, define $\bar{\chi }(\mathcal{O})\equiv \bar{\chi }(\mathcal{O}\cap \mathrm{D}_\mathrm{U})$, where $\bar{\chi }$ is defined in Eq. (\ref{Def chi bar}) on $\mathrm{Dom}(\bar{\chi })\equiv \mathrm{D}_\mathrm{U}$. Let $\mathcal{T}'$ be the set of the subsets $\mathcal{O} \subset \mathrm{N}_v$ such that 
\be\label{def-Topo-prime}
\forall \chi \in \mathcal{F}_v,\quad \bar{\chi }(\mathcal{O}) \mathrm{\ is\ an\ open\ set\ in\ }\mathbb{R}^N.
\ee
The set $\mathcal{T}'$ is a topology on $\mathrm{N}_v$.} \\

\noi {\it Proof.} This is an adaptation of the proof of Proposition C in Ref. \cite{A44}, replacing $\mathrm{M}$ by $\mathrm{N}_v$, $\mathrm{F}$ by $\mathcal{F}_v$, $\tilde{\chi }$ by $\bar{\chi }$, and $\mathbb{R}^3$ by $\mathbb{R}^N$. In particular, the proof that the whole set $\mathrm{N}_v$ (instead of $\mathrm{M}$) belongs to $\mathcal{T}'$ is exactly identical. Also, by definition of a nice $v$-adapted chart, the mapping $\bar{\chi}: \mathrm{D}_\mathrm{U} \rightarrow \mathbb{R}^N$ is injective. Therefore, if $\mathcal{O} _1\in \mathcal{T}'$ and $\mathcal{O} _2\in \mathcal{T}'$, we have
\bea
\bar{\chi } (\mathcal{O} _1\cap \mathcal{O} _2)& \equiv & \bar{\chi } ((\mathcal{O} _1\cap \mathcal{O} _2)\cap \mathrm{D}_\mathrm{U}) =  \bar{\chi } ((\mathcal{O} _1\cap\mathrm{D}_\mathrm{U})\cap ( \mathcal{O} _2\cap \mathrm{D}_\mathrm{U}))\nonumber\\
& = & \bar{\chi } (\mathcal{O} _1\cap\mathrm{D}_\mathrm{U})\cap \bar{\chi } ( \mathcal{O} _2\cap \mathrm{D}_\mathrm{U})\equiv \bar{\chi } (\mathcal{O} _1)\cap \bar{\chi } ( \mathcal{O} _2),
\eea
which is thus an open set of $\mathbb{R}^N$, for any $\chi \in \mathcal{F}_v$, so that $\mathcal{O} _1\cap \mathcal{O} _2\in \mathcal{T}'$. It is also trivial to check that the union of any family of subsets $\mathcal{O} _i\in \mathcal{T}'$ is still an element of $\mathcal{T}'$. \hfill $\square$\\

\noi To prove that the mappings $\bar{\chi}$ are continuous for this topology, and to prove the compatibility of any two mappings $\bar{\chi}$, $\bar{\chi}'$ on $\mathrm{N}_v$, associated with two nice $v$-adapted charts $\chi,\chi' \in \mathcal{F}_v$, the following difficulty arises: $\chi$ and $\chi'$ have in general different domains $\mathrm{U}$ and $\mathrm{U}'$. We may have $\mathrm{U} \cap \mathrm{U}'= \emptyset$, although there is some $l\in \mathrm{N}_v$ with
\be
 l \cap \mathrm{U} \ne \emptyset,\quad l\cap \mathrm{U}'\ne \emptyset.
\ee
I.e., it may happen that the domains of the charts $\chi $  and $\chi ' $ do not overlap, and that the domains of the mappings $\bar{\chi}$ and $ \bar{\chi}'$ do. To overcome this difficulty, we use the flow $F$ of the vector field $v$ to associate smoothly with any point $Y$ in some neighborhood $\mathrm{W}\subset \mathrm{U}$ of a point $X$,  a point $g(Y)\in \mathrm{U}' $:

\paragraph{Lemma.}\label{Lemma} {\it Let $v$ be a $\mathcal{C}^\infty $ vector field on $\mathrm{V}$. Let $\chi,\chi'\in \mathcal{F}_v$, with domains $\mathrm{U}$ and $\mathrm{U}'$, be such that $\mathrm{D}_\mathrm{U}\cap \mathrm{D}_{\mathrm{U}'} \ne \emptyset$. Let $l\in \mathrm{D}_\mathrm{U}\cap \mathrm{D}_{\mathrm{U}'}$ and $X\in l\cap \mathrm{U}$ and set $\chi (X)=(t,{\bf x})$, so that $\bar{\chi }(l)={\bf x}$. There is an open neighborhood $\Omega $ of ${\bf x}$ in $\mathbb{R}^N$ and a $\mathcal{C}^\infty $ mapping $g$ defined on an open neighborhood $\mathrm{W}\subset \mathrm{U}$ of $X$, such that for any ${\bf y}\in \Omega $ we have $(t,{\bf y})\in \chi (\mathrm{W})$ (so that ${\bf y}\in \bar{\chi } (\mathrm{D}_\mathrm{U})$), $\bar{\chi }^{-1}({\bf y})\in  \mathrm{D}_{\mathrm{U}'}$, and}
\be\label{chibar' rond chibar^-1}
 \forall {\bf y}\in \Omega,\quad \bar{\chi }'\circ \bar{\chi }^{-1}({\bf y})=P_S(\chi '(g(\chi ^{-1}(t,{\bf y})))).
\ee
\vspace{1mm}

\noi {\it Proof.}  Since $l\in \mathrm{D}_{\mathrm{U}'}$, there is some point $X'\in l\cap \mathrm{U}'$, and since $X\in l$ there is some $s\in \mathrm{I}_X$ such that $X'=F(s,X)$. Thus, the domain $\mathcal{D}$ of $F$ being open in $\mathbb{R}\times \mathrm{V}$ and $F$ being continuous, there is an interval $\mathrm{I}$ centered at $s $ and an open neighborhood $\mathrm{W}\subset \mathrm{U}$ of $X$ in $\mathrm{V}$, such that $\mathrm{I}\times\mathrm{W}\subset \mathcal{D}$ and $F(\mathrm{I}\times\mathrm{W})\subset \mathrm{U}'$. For $Y\in \mathrm{W}$, set $g(Y)=F(s,Y)$. This defines a  $\mathcal{C}^\infty $-mapping $g: \mathrm{W}\rightarrow \mathrm{U}' $. Because $\chi (\mathrm{W})$ is open in $\mathbb{R}^{N+1}$ and $(t,{\bf x})\in \chi (\mathrm{W})$, there is an interval $\mathrm{J}$ centered at $t $ and an open neighborhood $\Omega $ of ${\bf x}$ in $\mathbb{R}^N$, with $\mathrm{J}\times\Omega \subset \chi (\mathrm{W}) $. Hence, if ${\bf y}\in \Omega $, then indeed $(t,{\bf y})\in \chi (\mathrm{W})$, thus $Y\equiv \chi ^{-1}(t,{\bf y})\in \mathrm{W}\subset \mathrm{U}$, which implies that $l_Y\in \mathrm{D}_{\mathrm{U}}$ and that
\be
{\bf y} = P_S(\chi (Y))=\bar{\chi }(l_Y),
\ee
so ${\bf y}\in \bar{\chi }(\mathrm{D}_{\mathrm{U}})$. Moreover, we have $Y'\equiv g(Y)=F(s,Y)\in l_Y\cap \mathrm{U}'$, so $ l_Y =\bar{\chi }^{-1}({\bf y})\in \mathrm{D}_{\mathrm{U}'}$ and
\be
{\bf y}' \equiv P_S(\chi' (Y'))=\bar{\chi }' (l_Y),
\ee
whence follows (\ref{chibar' rond chibar^-1}). \hfill $\square$

\paragraph{Theorem 6.}\label{Theorem 6} {\it Let $v$ be a normal vector field on $\mathrm{V}$. Let $\mathcal{A}$ be the set of all mappings $\bar{\chi}$, where $\chi \in \mathcal{F}_v$. This set $\mathcal{A}$ is an atlas on the topological space $(\mathrm{N}_v,\mathcal{T}')$.}\\

\noi {\it Proof.} (i) Consider any $\chi '\in \mathcal{F}_v$, with domain $\mathrm{U}'$. Let us prove that $\bar{\chi }'$, which is defined on $\mathrm{D}_{\mathrm{U}'}$, is continuous for the topology induced on $\mathrm{D}_{\mathrm{U}'}$ by the topology $\mathcal{T}'$ on $\mathrm{N}_v$. Thus, $\mathrm{A}$ being any open set in $\mathbb{R}^N$, we must show that the set $\mathcal{O}_1\equiv \bar{\chi }'^{-1}(\mathrm{A})$ has the form $\mathcal{O}_1=\mathcal{O}\cap \mathrm{D}_{\mathrm{U}'}$, where $\mathcal{O}$ is such that we have (\ref{def-Topo-prime}). We shall actually show that $\mathcal{O}_1 \in \mathcal{T}'$, i.e., that we have (\ref{def-Topo-prime}) with $\mathcal{O}\equiv \mathcal{O}_1\subset \mathrm{D}_{\mathrm{U}'}$. To prove this, we may assume that $\mathcal{O}_1\ne \emptyset$. Moreover, let $\chi \in \mathcal{F}_v$, with domain $\mathrm{U}$. We may also assume that $\bar{\chi }(\mathcal{O}_1)\ne \emptyset$, i.e. $\mathcal{O}_1\cap \mathrm{D}_{\mathrm{U}} \ne \emptyset$, so let ${\bf x}\in \bar{\chi }(\mathcal{O}_1)$. We have to find a neighborhood $\mathrm{B} $ of ${\bf x}$ in $\mathbb{R}^N$, such that $\mathrm{B} \subset \bar{\chi }(\mathcal{O}_1)$, i.e., $\mathrm{B} \subset \bar{\chi }(\bar{\chi }'^{-1}(\mathrm{A}))$. Since ${\bf x}\in \bar{\chi }(\mathcal{O}_1)\equiv \bar{\chi }(\mathcal{O}_1\cap \mathrm{D}_{\mathrm{U}})$, we have that $l\equiv \bar{\chi }^{-1}({\bf x})\in (\bar{\chi }'^{-1}(\mathrm{A}))\cap \mathrm{D}_{\mathrm{U}}\subset \mathrm{D}_{\mathrm{U}'}\cap \mathrm{D}_{\mathrm{U}}$. In particular, there is some $X\in l\cap \mathrm{U}$, with $\chi (X)=(t,{\bf x})$ for some $t\in \mathbb{R}$. Hence, we may apply the \hyperref[Lemma]{Lemma} and get the corresponding open neighborhood $\Omega $ of ${\bf x}$. Thus (\ref{chibar' rond chibar^-1}) shows that the mapping ${\bf f}\equiv \bar{\chi }'\circ \bar{\chi }^{-1}$ is well defined and continuous over $\Omega $. Hence, since ${\bf x}\in \Omega $ and ${\bf x}'\equiv {\bf f}({\bf x})= \bar{\chi }'(l)\in \mathrm{A}$, which is open in $\mathbb{R}^N$, there is a neighborhood $\mathrm{B} \subset \Omega $ of ${\bf x}$ in $\mathbb{R}^N$, such that ${\bf f}(\mathrm{B})\subset \mathrm{A}$. Set ${\bf f}'\equiv \bar{\chi }\circ \bar{\chi }'^{-1}={\bf f}^{-1}$. For ${\bf y}\in \mathrm{B}$, we have thus ${\bf y}'\equiv {\bf f}({\bf y})\in \mathrm{A}$, that is, ${\bf y}={\bf f}'({\bf y}')\in {\bf f}'(\mathrm{A})=(\bar{\chi }\circ \bar{\chi }'^{-1})(\mathrm{A})$. Thus, for any  open set $\mathrm{A}$ in $\mathbb{R}^N$, the set $\mathcal{O}_1\equiv \bar{\chi }'^{-1}(\mathrm{A})$ belongs to $\mathcal{T}'$, as announced. This proves that $\bar{\chi }'$ is continuous for the topology induced on $\mathrm{D}_{\mathrm{U}'}$ by the topology $\mathcal{T}'$ on $\mathrm{N}_v$. Moreover, taking  $\mathrm{A}=\mathbb{R}^N$, we get that $\mathrm{D}_{\mathrm{U}'}=\bar{\chi }'^{-1}(\mathbb{R}^N)$ is open in $\mathrm{N}_v$, $\mathrm{D}_{\mathrm{U}'}\in \mathcal{T}'$. \\

\noi (ii) Given any $\chi \in \mathcal{F}_v$, with domain $\mathrm{U}$, let us show that the mapping $\bar{\chi }^{-1}:\bar{\chi }(\mathrm{D}_{\mathrm{U}})\rightarrow \mathrm{D}_{\mathrm{U}} \subset \mathrm{N}_v$, is continuous. Since $\mathrm{D}_{\mathrm{U}}$ is open as seen at the end of (i), we have to show that, for any $\mathcal{O}\in \mathcal{T}'$ such that $\mathcal{O}\subset \mathrm{D}_{\mathrm{U}}$, the set $\Omega \equiv (\bar{\chi }^{-1})^{-1}(\mathcal{O})$ is open in $\mathbb{R}^N$. But since $\mathcal{O}\subset \mathrm{D}_{\mathrm{U}}=\mathrm{Dom}(\bar{\chi })$, and since $\bar{\chi }$ is injective, we have $\Omega= \bar{\chi }(\mathcal{O})$. The fact that this is open in $\mathbb{R}^N$ follows from
the very definition of the topology $\mathcal{T}'$ in Eq. (\ref{def-Topo-prime}). With (i), this means that, for any $\chi \in \mathcal{F}_v$, the mapping $\bar{\chi }: \mathrm{D}_{\mathrm{U}} \rightarrow \bar{\chi }(\mathrm{D}_{\mathrm{U}})\subset \mathbb{R}^N$ is bicontinuous, thus is indeed a {\it chart} on the topological space $(\mathrm{N}_v,\mathcal{T}')$.\\

\noi (iii) Let us show that the domains of definition of the charts $\bar{\chi}$, for $\chi \in \mathcal{F}_v$, cover the whole set $\mathrm{N}_v$. Given $l\in \mathrm{N}_v$, let $X\in l$. Since $v$ is a normal vector field on $\mathrm{V}$, we know from \hyperref [Theorem 5]{Theorem 5} that there is a nice $v$-adapted chart $\chi $ whose domain $\mathrm{U}$ is a neighborhood of $X$. Thus $X\in l\cap \mathrm{U}$, so $l\in \mathrm{D}_{\mathrm{U}}$, Q.E.D.\\

\noi (iv) Finally, given any two nice $v$-adapted charts $\chi,\chi' \in \mathcal{F}_v$, having domains $\mathrm{U}$ and $\mathrm{U}'$ respectively, let us show the compatibility of the two charts $\bar{\chi}$, $\bar{\chi}'$ on $\mathrm{N}_v$, with domains $\mathrm{D}_{\mathrm{U}}$ and $\mathrm{D}_{\mathrm{U}'}$. The relevant case is when $\mathrm{D}_{\mathrm{U}}\cap \mathrm{D}_{\mathrm{U}'} \ne \emptyset$, so that $\mathrm{Dom}( \bar{\chi }'\circ \bar{\chi }^{-1})=\bar{\chi }(\mathrm{D}_{\mathrm{U}}\cap \mathrm{D}_{\mathrm{U}'} )\ne \emptyset$. Thus we may apply the \hyperref[Lemma]{Lemma}. Its Eq. (\ref{chibar' rond chibar^-1}) shows that, given any ${\bf x} \in \bar{\chi }(\mathrm{D}_{\mathrm{U}}\cap \mathrm{D}_{\mathrm{U}'} )$, it has a neighborhood $\Omega $ in which the function $\bar{\chi }'\circ \bar{\chi }^{-1}$ is given as a composition of $\mathcal{C}^\infty $ functions, hence $\bar{\chi }'\circ \bar{\chi }^{-1}$ is $\mathcal{C}^\infty $ on its domain. \hfill $\square$\\

\vspace{2mm}
\noi Thus, we have endowed the set $\mathrm{N}_v$ with first the topology $\mathcal{T}'$ defined by (\ref{def-Topo-prime}), and then with a canonical atlas $\mathcal{A}$ of compatible charts, which are simply the mappings $\bar{\chi}$, where $\chi \in \mathcal{F}_v$ is any nice $v$-adapted chart. To call this a differentiable manifold in the rather usual sense of Note \ref{Def Dif Manifold} needs that the topological space $(\mathrm{N}_v,\mathcal{T}')$ be metrizable and separable --- hence, in particular, that it be Hausdorff. We do not have very general results on the latter point.

\paragraph{Proposition 7.}\label{Proposition 7} {\it  Let $v$ be a normal vector field on $\mathrm{V}$.} (i) {\it Any two points $l \ne l'$ in the orbit space $\mathrm{N}_v$ are topologically distinguishable, i.e., there exists an open set $\mathcal{O}\in \mathcal{T}'$ such that $l\in \mathcal{O}$ and $l' \notin \mathcal{O}$.}\\
\noi (ii) {\it Suppose $l \ne l'$ but $l$ and $l'$ both belong to the domain $\mathrm{D}_{\mathrm{U}}$ of a chart $\bar{\chi } $, where $\chi \in \mathcal{F}_v$ with $\mathrm{Dom}\,\chi =\mathrm{U}$. Then $l$ and $l'$ are separated by neighborhoods, i.e., there are two open sets $\mathcal{O},\mathcal{O}'\in \mathcal{T}'$, such that $l\in \mathcal{O}$, $l'\in \mathcal{O}'$, and $\mathcal{O}\cap \mathcal{O}'=\emptyset$.} \\
\noi (iii) {\it Suppose there is a chart $\chi \in \mathcal{F}_v$, such that any maximal integral curve $l\in \mathrm{N}_v$ intersects its domain $\mathrm{U}$, so that $\mathrm{D}_{\mathrm{U}}=\mathrm{N}_v$. Then the topological space $(\mathrm{N}_v,\mathcal{T}')$ is metrizable and separable.} [Hence, in particular, it has the Hausdorff property, i.e., any two distinct elements $l,l'\in \mathrm{N}_v$ are separated by neighborhoods --- as follows also from Point (ii).]\\

\noi {\it Proof.} (i) Let $X\in l$. There is some open neighborhood $\mathrm{U}$ of $X$, such that $\mathrm{U}\cap l'=\emptyset$: if that were not true, then, taking any distance on $\mathrm{V}$ that defines its topology, and considering $\mathrm{U}_n\equiv \mathrm{B}(X,1/n)$, we would get a sequence $(X_n)$ with $X_n\in l'$ and $X_n\rightarrow X$, hence $X\in l'$ since we defined that a normal vector field has all its maximal integral curves closed; but this implies that $l=l'$, which is a contradiction. By \hyperref[Theorem 5]{Theorem 5}, let $\chi \in \mathcal{F}_v $ whose domain $\mathrm{U}_2$ is an open neighborhood of $X$, and set $\mathrm{U}_1\equiv \mathrm{U}\cap \mathrm{U}_2$. The restriction $\chi _1$ of $\chi $ to $\mathrm{U}_1\subset \mathrm{U}$ is still a nice $v$-adapted chart: it is obviously $v$-adapted, we have $\mathrm{D}_{\mathrm{U}_1}\subset \mathrm{D}_{\mathrm{U}}$, and the mapping $\bar{\chi }_1: \mathrm{D}_{\mathrm{U}_1}\rightarrow \mathbb{R}^N$ is clearly the restriction of $\bar{\chi }$ to $ \mathrm{D}_{\mathrm{U}}$, hence it also is injective. Therefore, $\mathcal{O}\equiv \mathrm{D}_{\mathrm{U}_1}$ is an open set, $\mathcal{O}\in \mathcal{T}'$. Since $X\in l\cap \mathrm{U}_1$, we have $l\in \mathcal{O}$; and since $l'\cap \mathrm{U}_1\subset l'\cap \mathrm{U}= \emptyset$, we have $l' \notin \mathcal{O}$.\\

\noi (ii) Since $l\in \mathrm{D}_{\mathrm{U}}$ and $l'\in \mathrm{D}_{\mathrm{U}}$ with $l\ne l'$, and since $\bar{\chi } $ is defined on $\mathrm{D}_{\mathrm{U}}$ and injective, we have
\be
{\bf x}\equiv \bar{\chi } (l) \ne {\bf x}'\equiv \bar{\chi } (l').
\ee
Let $\Omega $ and $\Omega '$ be open neighborhoods in $\mathbb{R}^N$ of ${\bf x}$ and ${\bf x}'$ respectively, such that $\Omega \cap \Omega '= \emptyset$. Set $\mathcal{O}\equiv \bar{\chi }^{-1}(\Omega )$ and $\mathcal{O}'\equiv \bar{\chi }^{-1}(\Omega' )$. These are open sets such that $l\in \mathcal{O}$ and $l'\in \mathcal{O} '$, and we have $ \mathcal{O}\cap \mathcal{O}'=\bar{\chi }^{-1}(\Omega \cap \Omega ')= \emptyset$.\\

\noi (iii) This follows from the fact that $\bar{\chi }$ is a homeomorphism of its domain $\mathrm{D}_{\mathrm{U}}$ onto its range $\bar{\chi }(\mathrm{D}_{\mathrm{U}})\subset \mathbb{R}^N$.  \hfill $\square$\\

\noi {\bf Example.} The assumption made for Point (iii) is fulfilled, in particular, in the following case, which occurs frequently in relativistic theories of gravitation. Assume a chart $(\chi,\mathrm{U}) $ is defined on the whole of the manifold: $\mathrm{U}=\mathrm{V}$, which means that $\chi $ is a diffeomorphism of $\mathrm{V}$ onto the open subset $\Gamma \equiv \chi (\mathrm{V})$ of $\mathbb{R}^{N+1}$. Then, the tangent vector field $v$ to the world lines $l_{{\bf a}}$ given by (\ref{world line}) for $N=3$, with constant component vector ${\bf v}_0=(1,0,...,0)$ in the chart $\chi$, is the pushforward vector field of the constant vector field ${\bf v}({\bf X})={\bf v}_0$ for ${\bf X} \in \Gamma $ by the diffeomorphism $\chi ^{-1}$. Hence, by No. (iv) in the \hyperref[Examples normal]{Examples} above, $v$ is a normal vector field on $\mathrm{V}$. Due to No. (ii) in these examples, the orbits of $v$ are the images of the orbits of ${\bf v}$ by $\chi ^{-1}$, hence [by No. (iii)] are the connected components of the lines $l_{{\bf a}}\ [{\bf a}=(a^j)\in P_S(\Gamma )]$. Hence, the chart $\chi $ is $v$-adapted, for $P_S(\chi (X))={\bf a} $ if $X\in l_{{\bf a}}$. If actually all lines $l_{{\bf a}}$ defined in (\ref{world line}) are {\it connected} (which happens iff the domain of the time coordinate $x^0$ is an interval for any such line), then $\chi $ is nice (${\bf a}={\bf a}'\Rightarrow l_{{\bf a}}=l_{{\bf a}'}$). Thus, in that case, $\chi \in \mathcal{F}_v$. Since the domain of $\chi $ is $\mathrm{U}=\mathrm{V}$, we have then $\mathrm{D}_{\mathrm{U}}= \mathrm{N}_v$, so $\mathrm{N}_v$ is metrizable and separable. This case includes of course standard situations, e.g. an inertial frame (e.g. with Cartesian coordinates) in Minkowski spacetime; a uniformly rotating frame (e.g. with ``rotating Cartesian coordinates" \cite{A47}) in Minkowski spacetime [even though in that case the lines (\ref{world line}) are spacelike when $\rho \equiv \sqrt{x^2+y^2}> c/\omega$]; harmonic coordinates in an asymptotically flat spacetime \cite{Fock1959}; etc. It also includes known singular solutions of general relativity such as the singular Schwarzschild-Kruskal-Szekeres spacetime: the Kruskal-Szekeres coordinates  $(T,\xi ,\theta ,\phi )$ \cite{Kruskal1960,Szekeres1960} cover the whole of the ``maximally extended" Schwarzschild manifold. Since the domain of the coordinates $T,\xi $ is: $\xi \in \mathbb{R}, T^2-\xi ^2<1$, i.e. $T\in ]-\sqrt{1+\xi ^2},+\sqrt{1+\xi ^2}[$, each line $l_{{\bf a}}$ [with ${\bf a}\equiv (\xi ,\theta ,\phi )$] is connected. Thus this global chart on the Schwarzschild spacetime does define a global space manifold. Moreover, the tangent vector field $v$ to these lines (\ref{world line}) is time-like. 

\paragraph{Proposition 8.}\label{Proposition 8} {\it Assume that $v$ is a normal vector field on $\mathrm{V}$.} (i) {\it There is a countable cover of $\mathrm{V}$ by open sets $\mathrm{U}_n$ such that, for any integer $n$, there is a nice $v$-adapted chart, $\chi _n \in \mathcal{F}_v$, having domain $\mathrm{U}_n$.} (ii) {\it Then, setting $\mathrm{D}_n\equiv \mathrm{D}_{\mathrm{U}_n}$, the sequence $(\mathrm{D}_n)$ is a countable cover of $\mathrm{N}_v$ by metrizable open subsets. Hence the topological space $(\mathrm{N}_v, \mathcal{T}')$ is separable.}\\

\noi {\it Proof.} (i) By \hyperref[Theorem 5]{Theorem 5}, for any $X\in \mathrm{V}$ there is a nice $v$-adapted chart $\chi_X \in \mathcal{F}_v $, such that its domain $\mathrm{U}_X$ is an open neighborhood of $X$. 
But, since  $\mathrm{V}$ is metrizable and separable, there exists a countable basis $(\mathrm{V}_n)_{n\in \mathbb{N}}$ for the open sets of $\mathrm{V}$. Hence, for any $X\in \mathrm{V}$, there is some integer $\widetilde{n}(X)$ such that 
\be\label{V_n & U_X}
X\in \mathrm{V}_{\widetilde{n}(X)}\subset \mathrm{U}_X. 
\ee
This defines a mapping $\widetilde{n}: \mathrm{V}\rightarrow \mathbb{N}$ 
and we have
\be\label{V union V_n}
\mathrm{V}=\bigcup _{n\in \widetilde{n}(\mathrm{V})} \mathrm{V}_n.
\ee
We may define a mapping $\widetilde{n}(\mathrm{V})\rightarrow \mathrm{V}, \ n\mapsto X_n$, by choosing $X_n$, for any $n \in \widetilde{n}(\mathrm{V})$, as being one of the points $X\in \mathrm{V}$ such that $n= \widetilde{n}(X)$. From (\ref{V_n & U_X}), it follows then that, for any $n \in \widetilde{n}(\mathrm{V})$, we have
\be\label{V_n subset U_n}
\mathrm{V}_n=\mathrm{V}_{\widetilde{n}(X_n)}\subset \mathrm{U}_{X_n}.
\ee
For $n \in \widetilde{n}(\mathrm{V})$, define $\mathrm{U}_n\equiv \mathrm{U}_{X_n}$ and $\chi _n\equiv \chi _{X_n}\in \mathcal{F}_v$. From (\ref{V union V_n}) and (\ref{V_n subset U_n}), it results that the countable family $(\mathrm{U}_n)_{n \in \widetilde{n}(\mathrm{V})}$ is as in Statement (i).\\

(ii) Note first that, since $\chi _n \in \mathcal{F}_v$, with domain $\mathrm{U}_n$, it follows from \hyperref[Theorem 6]{Theorem 6} that $\mathrm{D}_n$, the domain of the associated chart $\bar{\chi }_n$ on $\mathrm{N}_v$, is open in $\mathrm{N}_v$. If $l\in \mathrm{N}_v$, let $X\in l$ and, since $(\mathrm{U}_n)$ is a cover of $\mathrm{V}$, let $n$ be such that $X\in \mathrm{U}_n$. We have thus $l\cap \mathrm{U}_n \ne \emptyset$, i.e. $l\in \mathrm{D}_n$. So $(\mathrm{D}_n)$ is a countable open cover of $\mathrm{N}_v$. Since  $\bar{\chi }_n$ is a homeomorphism of $\mathrm{D}_n$ onto $\bar{\chi }_n(\mathrm{D}_n)\subset \mathbb{R}^N$, it follows that $\mathrm{D}_n$ is a metrizable and separable space. Hence it is second-countable, i.e., there exists a countable basis $(\mathcal{O}_{nm})_{m\in \mathbb{N}}$ for the open subsets of $\mathrm{D}_n$. Since any open subset $\mathcal{O}$ of $\mathrm{N}_v$ is the countable union of the open subsets $\mathcal{O}_n\equiv \mathcal{O}\cap \mathrm{D}_n$ of $\mathrm{D}_n$, we have that $(\mathcal{O}_{nm})_{n,m\in \mathbb{N}}$ is a countable basis for the topology $\mathcal{T}'$ of $\mathrm{N}_v$. Thus also $\mathrm{N}_v$ is second-countable, hence it is separable.\hfill $\square$

\section{The local manifold as an open subset of the global one}\label{Local_vs_global}

Let $v$ be a \hyperref[Normal]{normal vector field} on $\mathrm{V}$. In addition, as in Subsect. \ref{Local}, let $\mathrm{F}$ be a (local) reference frame, thus an equivalence class of charts for the relation (\ref{R_U}), in which $\mathrm{U}$ is a given open subset of $\mathrm{V}$. 
\footnote{\
As in Sects. \ref{AdaptedCharts} and \ref{N_v manifold}, the dimension of $\mathrm{V}$ is $N+1$, where $N$ is any integer $\ge 1$. All results summarized in Subsect. \ref{Local} hold true if one substitutes any integer $N \ge 1$ for the integer $3$, and $N+1$ for $4$ \cite{A44}.
} 
In Subsect. \ref{Local}, the local space manifold $\mathrm{M}_\mathrm{F}$ associated with $\mathrm{F}$ was defined as the set of the world lines (\ref{l-in-M-by-P_S}). On the other hand, the orbit set $\mathrm{N}_v$ defined in Subsect. \ref{Definitions}: the set of the maximal integral curves of $v$, was endowed in Sect. \ref {N_v manifold} with a topology $\mathcal{T}'$ and an atlas $\mathcal{A}$, which (assuming that  $\mathcal{T}'$ is metrizable and separable) makes it a differentiable manifold. Thus, we have also a global space manifold: $\mathrm{N}_v$. When the charts $\chi \in \mathrm{F}$, all having domain $\mathrm{U}$, are nice $\,v$--adapted charts, i.e. belong to  $\mathcal{F}_v$, we have the following tight relation between $\mathrm{M}_\mathrm{F}$ and $\mathrm{N}_v$:

\paragraph{Theorem 7.}\label{Theorem 7} {\it Assume that $\mathrm{F}\subset \mathcal{F}_v$. For any $l\in \mathrm{M}_\mathrm{F}$, there is a unique maximal integral curve $l'\in \mathrm{N}_v$ such that, for any $X\in l$, we have $l'=l_X$. It holds $l=l'\cap \mathrm{U}$. The mapping $I:\, l\mapsto l'$ is a diffeomorphism of $\mathrm{M}_\mathrm{F}$ onto the open subset $\mathrm{D}_\mathrm{U}$ of $\mathrm{N}_v$.}\\

\noi {\it Proof.} Let $l\in \mathrm{M}_\mathrm{F}$. By the definition of $\mathrm{M}_\mathrm{F}$ near Eq. (\ref{l-in-M-by-P_S}), there is some chart $\chi\in \mathrm{F}$ and some ${\bf x}\in P_S(\chi(\mathrm{U}))\subset \mathbb{R}^N$, such that 
\be\label{def l}
l =\{\,X\in \mathrm{U};\ P_S(\chi (X))={\bf x}\,\}.
\ee
Let $\,X_1 \in l$ and $\,X_2 \in l$. Denote the maximal integral curves of $v$ 
at $X_1$ and $X_2$ as $l'_1\equiv l_{X_1}$, $l'_2\equiv l_{X_2}$. Since $\chi$ 
is $v$-adapted, there exist ${\bf x}_1, {\bf x}_2\in \mathbb{R}^N$, such that 
$\forall X\in l'_1\cap \mathrm{U},\ P_S(\chi (X))= {\bf x}_1$; and $\forall X\in l'_2\cap \mathrm{U},\ P_S(\chi (X))= {\bf x}_2$. In particular, since $X_j\in l'_j\cap \mathrm{U}$, we have $P_S(\chi (X_j))= {\bf x}_j \ (j=1,2)$. But since $X_j\in l$, we have also $P_S(\chi (X_j))={\bf x}$ by (\ref{def l}), so ${\bf x}={\bf x}_1={\bf x}_2$. Because the $v$-adapted chart $\chi$ is nice, it follows that $l'_1=l'_2$. Therefore, we define a mapping $I:\mathrm{M}_\mathrm{F}\rightarrow \mathrm{N}_v $ by associating with any $l\in \mathrm{M}_\mathrm{F}$, the unique maximal integral curve $l'\in \mathrm{N}_v$, such that for any $X\in l$, we have $l'=l_X$. Note that actually $l'\in \mathrm{D}_\mathrm{U}$. Owing to the definitions (\ref{def-chi-tilde}) and (\ref{Def chi bar}), we have ${\bf x}= \widetilde{\chi }(l)=\bar{\chi}(l')$, and since here $l$ is any element of $\mathrm{M}_\mathrm{F}$, this shows that $\widetilde{\chi }(\mathrm{M}_\mathrm{F})\subset \bar{\chi }(\mathrm{D}_\mathrm{U})$. Thus
\be\label{I=chi bar ^-1 rond chi tilde}
I(l)=\bar{\chi }^{-1}(\widetilde{\chi }(l)),
\ee
so that we have simply $I=\bar{\chi }^{-1}\circ \widetilde{\chi }$, for whatever chart $\chi \in \mathrm{F}$. Let us show that $\,l=l' \cap\mathrm{U}$. By definition, for any $X\in l$, we have $l'=l_X$, hence $X\in l'$, and since $l\subset \mathrm{U}$ by the definition (\ref{def l}), we have $\,l\subset l' \cap\mathrm{U}$. Conversely, consider any $\chi \in \mathrm{F}$; this is by assumption a $v$-adapted chart and, as we showed before (\ref{I=chi bar ^-1 rond chi tilde}), we have $\bar{\chi}(l')=\widetilde{\chi }(l)\equiv {\bf x}$. Therefore, by the definition (\ref{Def adapted}), we have for any $X\in l'\cap \mathrm{U}$: $P_S(\chi (X))={\bf x}$. Then (\ref{def l}) implies that $X\in l$, so $l' \cap\mathrm{U}\subset l$. \\

As we showed, the mapping $I$ is defined on the whole of $\mathrm{M}_\mathrm{F}$ and ranges into $\mathrm{D}_\mathrm{U}$, which is an open subset of $\mathrm{N}_v$. Let us show that in fact $I(\mathrm{M}_\mathrm{F})=\mathrm{D}_\mathrm{U}$. Let $l'\in \mathrm{D}_\mathrm{U}$, so there exists $X\in l'\cap \mathrm{U}$. Let $\chi \in \mathrm{F}$ and set ${\bf x}\equiv P_S(\chi (X))$ and $l\equiv \{\,Y\in \mathrm{U};\ P_S(\chi (Y))={\bf x}\,\}$. Clearly $l\in \mathrm{M}_\mathrm{F}$ and $X\in l$. By the definition of $I$, we have that, for any $Y\in l$, $l_Y=I(l)$. In particular, $l_X=I(l)$. But since $X\in l'$, we have $l_X=l'$, hence $l'=I(l)$: thus indeed $\mathrm{D}_\mathrm{U}\subset I(\mathrm{M}_\mathrm{F})$. Note that again here, from the definitions (\ref{def-chi-tilde}) and (\ref{Def chi bar}) we have ${\bf x}= \widetilde{\chi }(l)=\bar{\chi}(l')$, and since now $l'$ is any element of $\mathrm{D}_\mathrm{U}$, this shows that $\bar{\chi }(\mathrm{D}_\mathrm{U})\subset \widetilde{\chi }(\mathrm{M}_\mathrm{F})$. Since the reverse inclusion has been proved before (\ref{I=chi bar ^-1 rond chi tilde}), we have $\bar{\chi }(\mathrm{D}_\mathrm{U})= \widetilde{\chi }(\mathrm{M}_\mathrm{F})$.\\

As shown in Ref. \cite{A44}, $\widetilde{\chi }$ is a global chart on the differentiable manifold $\mathrm{M}_\mathrm{F}$, for any $\chi \in \mathrm{F}$. As shown in \hyperref[Theorem 6]{Theorem 6}, $\bar{\chi }$ is a chart with domain $\mathrm{D}_\mathrm{U}$ on the differentiable manifold $\mathrm{N}_v$, also for any $\chi \in \mathrm{F}$. Moreover, as we just saw, we have $\bar{\chi }(\mathrm{D}_\mathrm{U})= \widetilde{\chi }(\mathrm{M}_\mathrm{F})$. It follows that the one-to-one mapping $I=\bar{\chi }^{-1}\circ \widetilde{\chi }$, from $\mathrm{Dom}(\widetilde{\chi })=\mathrm{M}_\mathrm{F}$ onto $I(\mathrm{M}_\mathrm{F})=\mathrm{D}_\mathrm{U}=\mathrm{Dom}(\bar{\chi })$, is a diffeomorphism. Therefore, $I$ is an immersion of $\mathrm{M}_\mathrm{F}$ into $\mathrm{N}_v$. Actually, recall that $\mathrm{D}_\mathrm{U}$ is more specifically an {\it open subset} of $\mathrm{N}_v$. \hfill $\square$\\

\vspace{1mm}
\noi Thus, if $\mathrm{F}\subset \mathcal{F}_v$, {\it we may identify the local space $\mathrm{M}_\mathrm{F}$ with the open subset  $I(\mathrm{M}_\mathrm{F})=\mathrm{D}_\mathrm{U}$ of the global space $\mathrm{N}_v$.} Since we proved that $\,l=I(l) \cap\mathrm{U}$, we can say that {\it $\mathrm{M}_\mathrm{F}$ is made of the intersections with the local domain $\mathrm{U}$ of the maximal integral curves of $v$.} Given that each world line $l\in \mathrm{M}_\mathrm{F}$ is invariant under any exchange of the chart $\chi \in \mathrm{F}$ for another chart $\chi' \in \mathrm{F}$, to say that $\mathrm{F}\subset \mathcal{F}_v$ is equivalent to say that {\it one} chart $\chi \in \mathrm{F}$ is a nice $v$-adapted chart.\\

\vspace{2mm}
\noi {\bf Acknowledgement.} Point (ii) in  \hyperref[Proposition 0]{Proposition 0}, the possibility of proving Point (ii) in \hyperref[Theorem 3]{Theorem 3} by noting that $F$ is a submersion (cf. \hyperref[Remark 3.1]{Remark 3.1}), as well as the proof of Point (i) in  \hyperref[Proposition 8]{Proposition 8} (rather than taking it as an additional assumption), were suggested by the referee while commenting on the first version. The new Point (ii) in  Proposition 0 led me to an improvement of \hyperref[Theorem 5]{Theorem 5}. Frank Reif\mbox{}ler suggested to check the example of a rotating frame.\\



\begin{thebibliography}{9}
\small


\bibitem{Lachièze-Rey2003}
M. Lachi\`eze-Rey, ``Space and spacetime," {\it Cosmology and gravitation: Xth Brazilian School of Cosmology and Gravitation} (AIP Conference Proceedings, Vol. 668, 2003), pp. 173-186.

\bibitem{Cattaneo1958}
C. Cattaneo, ``General relativity: relative standard mass, momentum, energy and gravitational field in a general system of reference," {\it il Nuovo Cimento} {\bf 10}, 318--337 (1958).

\bibitem{Massa1974a}
E. Massa, ``Space tensors in general relativity. I. Spatial tensor algebra and analysis," {\it Gen. Rel. Grav.} {\bf 5}, 555--572 (1974).

\bibitem{Massa1974b}
E. Massa, ``Space tensors in general relativity. II. Physical applications," {\it Gen. Rel. Grav.} {\bf 5}, 573--591 (1974).

\bibitem{Mitskievich1996}
N. V. Mitskievich, {\it Relativistic physics in
arbitrary reference frames} (Nova Science Publishers, Hauppauge, NY, 2007). 

\bibitem{JantzenCariniBini1992}
R. T. Jantzen, P. Carini, and D. Bini, ``The many faces of gravitoelectromagnetism," {\it Ann. Phys. (New York)} {\bf 215}, 1--50 (1992). 

 \bibitem{A16} M. Arminjon, ``On the extension of Newton's second law to theories of gravitation in curved space-time," {\it Arch. Mech.} {\bf 48}, 551--576 (1996). 

\bibitem{A44}
M. Arminjon and F. Reif\mbox{}ler, ``General reference frames and their associated space manifolds," {\it Int. J. Geom. Meth. Mod. Phys.} {\bf 8}, 155--165 (2011). 

\bibitem{A43}
M. Arminjon and F. Reif\mbox{}ler, ``A non-uniqueness problem of the Dirac theory in a curved spacetime," {\it Ann. Phys. (Berlin)} {\bf 523}, 531--551 (2011). 

\bibitem{A47}
M. Arminjon, ``A solution of the non-uniqueness problem of the Dirac Hamiltonian and energy operators," {\it Ann. Phys. (Berlin)} {\bf 523}, 1008--1028 (2011). 


\bibitem{DieudonneTome3}
J. Dieudonn\'e, {\it El\'ements d'Analyse,} Tome 3 (2nd French edition, Gauthier-Villars, Paris 1974). (English edition: {\it Treatise on Analysis,} Volume 3, Academic Press, New York, 1972.) 

\bibitem{DieudonneTome4}
J. Dieudonn\'e, {\it El\'ements d'Analyse,} Tome 4 (1st edition, Gauthier-Villars, Paris, 1971), pp. 4--7. (English edition: {\it Treatise on Analysis,} Volume 4, Academic Press, New York, 1974.) 

\bibitem{AbrahamMarsdenRatiu}
R. Abraham, J. E. Marsden and T. S. Ra\c tiu, {\it Manifolds, Tensor Analysis, and Applications} (2nd edition, Springer, 1988), p. 216.

\bibitem{DieudonneTome4-2}
Ref. \cite{DieudonneTome4}, Problem 7 on pp. 11-12.


\bibitem{DieudonneTome1} J. Dieudonn\'e, {\it El\'ements d'Analyse,} Tome 1 : {\it Fondements de l'Analyse Moderne} (2nd edition, Gauthier-Villars, Paris, 1972). (English edition: {\it Treatise on Analysis,} Volume 1: {\it Foundations of Modern Analysis,} Revised edition, Academic Press, New York, 1969.) 

\bibitem{GuilleminPollack-p28}
V. Guillemin, A. Pollack, {\it Differential Topology} (Prentice-Hall, Englewood Cliffs, NJ, 1974), Theorem on p. 28.

\bibitem{CalcutGompf2013}
J. S. Calcut and R. E. Gompf, ``Orbit spaces of gradient vector fields," {\it Ergodic Theory and Dynam. Syst.} {\bf 33}, 1732--1747 (2013).

\bibitem{Fock1959}
V. I. Fock, {\it The Theory of Space, Time and Gravitation} (First
English edition, Pergamon, Oxford, 1959). (Russian original edition:
 {\it Teoriya Prostranstva, Vremeni i Tyagotenie}, Gos. Izd.
Tekhn.-Teoret. Liter., Moskva, 1955.)

\bibitem{Kruskal1960}
M. D. Kruskal, ``Maximal extension of the Schwarzschild metric," {\it Phys. Rev.} {\bf 119}, 1743--1745 (1960).

\bibitem{Szekeres1960}
G. Szekeres, ``On the singularities of a Riemannian manifold," {\it  Publ. Math. Debrecen} {\bf 7}, 285--301 (1960). 




\end{thebibliography}
\end{document}